\documentclass[12pt]{article}
\usepackage{amsmath,amssymb,amsbsy}
\usepackage{graphicx}
\usepackage{multirow}
\usepackage{authblk}
\usepackage{algorithm}
\usepackage{xcolor}
\usepackage[normalem]{ulem}
\usepackage{tikz}
\usepackage{tkz-graph}
\usepackage[letterpaper, left=2.5cm, right=2.5cm, top=2.5cm,bottom=2.5cm,dvips]{geometry}

\newtheorem{proposition}{Proposition}[section]
\newtheorem{theorem}[proposition]{Theorem}

\newtheorem{corollary}[proposition]{Corollary}

\newtheorem{definition}[proposition]{Definition}

\newtheorem{remark}[proposition]{Remark}
\newtheorem{problem}[proposition]{Problem}

\newcommand{\qed}{\hfill \rule{.1in}{.1in}}

\makeatletter
\def\imod#1{\allowbreak\mkern10mu({\operator@font mod}\,\,#1)}
\makeatother

\begin{document}

\title{On b-acyclic chromatic number of a graph}

\author[1]{Marcin Anholcer\thanks{Partially supported by the National Science Center of Poland under grant no. 2020/37/B/ST1/03298.}}
\author[2]{Sylwia Cichacz\thanks{The work of the author was partially supported by the Faculty of Applied Mathematics AGH UST statutory tasks within subsidy of Ministry of Education and Science. The author would like to thank prof. Agnieszka Görlich for her valuable comments.}}
\author[3]{Iztok Peterin\thanks{The author is partially supported by the Slovenian Research Agency by the projects No. J1-1693 and J1-9109 and is also with the Institute of Mathematics, Physics and Mechanics, Jadranska 19, 1000 Ljubljana, Slovenia}}

\affil[1]{\scriptsize{}Pozna\'n University of Economics and Business, Institute of Informatics and Quantitative Economy}
\affil[ ]{Al.Niepodleg{\l}o\'sci 10, 61-875 Pozna\'n, Poland, \textit{m.anholcer@ue.poznan.pl}}
\affil[ ]{}
\affil[2]{AGH University of Science and Technology, Faculty of Applied Mathematics}
\affil[ ]{Al. Mickiewicza 30, 30-059 Krak\'ow, Poland, \textit{cichacz@agh.edu.pl}}
\affil[ ]{}
\affil[3]{University of Maribor, Faculty of Electrical Engineering and Computer Science}
\affil[ ]{Koro\v{s}ka 46, 2000 Maribor, Slovenia,\textit{iztok.peterin@um.si}}

\maketitle

\begin{abstract}
Let $G$ be a graph. We introduce the acyclic b-chromatic number of $G$ as an analogue to the b-chromatic number of $G$. An acyclic coloring of a graph $G$ is a map $c:V(G)\rightarrow \{1,\dots,k\}$ such that $c(u)\neq c(v)$ for any $uv\in E(G)$ and the induced subgraph on vertices of any two colors $i,j\in \{1,\dots,k\}$ induces a forest. On the set of all acyclic colorings of $G$ we define a relation whose transitive closure is a strict partial order. The minimum cardinality of its minimal element is then the acyclic chromatic number $A(G)$ of $G$ and the maximum cardinality of its minimal element is the acyclic b-chromatic number $A_b(G)$ of $G$. We present several properties of $A_b(G)$. In particular, we derive $A_b(G)$ for several known graph families, derive some bounds for $A_b(G)$, compare $A_b(G)$ with some other parameters and generalize some influential tools from b-colorings to acyclic b-colorings. 
\end{abstract}

  \medskip


\section{Introduction}\label{sec_Intro}

The computation of chromatic number $\chi(G)$ of a graph $G$ is a well know difficult problem that is NP-hard. As such one often tries to find some approximate values for it. For this one first needs a (proper) coloring of $G$. There are some simple approaches how to find such a coloring. Let us recall two of them. First is the greedy approach, also sometimes called first-fit, where one starts with totally uncolored graph. Vertices are then colored in some arbitrary order by the rule that an uncolored vertex receives the minimum color that is not present in its neighborhood until that moment. At the end of this procedure we obtain a coloring of $G$ and the number of colors is an upper bound for $\chi(G)$. The other approach starts at the other end where all vertices are colored by a different color from $\{1,\dots,|V(G)|\}$. In what follows we need to find at each step a color class that is without a vertex having neighbors of all the remaining colors. Every vertex of such a class can then be recolored with some color, the one that is missing in its neighborhood, and we obtain a new coloring with one color less than before. We stop with this when every color class has a vertex with all the other colors in its neighborhood. Again the number of the colors at the last stage is an upper bound for $\chi(G)$. 

Both mentioned procedures can result in a coloring with the number of colors that are close to $\chi(G)$ and if we are lucky, then even with $\chi(G)$ colors. However, the difference between the obtained number of colors and $\chi(G)$ can also be arbitrarily large. For both a wide range of studies deals with the worst case. In the greedy approach the worst number of colors that can be obtained is called the Grundy chromatic number $\Gamma (G)$ of $G$ and the worst case in the other presented procedure is called the b-chromatic number $\varphi (G)$ of $G$.

The Grundy chromatic number was introduced by Christen and Selkow \cite{ref_ChrSel} and then investigated by numerous authors. Let us cite just a few results. Erd\"os et al. \cite{ref_ErdHarHedLas} proved that for every finite graph the Grundy chromatic number is equal to the ochromatic number (the one corresponding to the parsimonious proper coloring). Telle and Proskurowski in \cite{ref_TelPro} presented the first polynomial-time algorithm for computing the Grundy number of partial $k$-trees. DeVilbiss, Johnson and Matzke \cite{ref_DeVJohMat} proved the values of this graph invariant for the line graphs of the regular Tur\'an graphs. The complexity of finding Grundy number was analyzed e.g. by  Zaker \cite{ref_Zak} and Bonnet et al. \cite{ref_BonFouKimSik}.

The b-chromatic number was introduced by Irving and Manlove in 1999 \cite{irma-99}. They have shown that determining the b-chromatic number of a graph is an NP-hard problem. The problem is still NP-hard for connected bipartite graphs as shown by Kratochv{\'{\i}}l, Tuza and Voigt \cite{krtuvo-02}. In contrast, the exact result for $\varphi (T)$ for every tree $T$ was presented already in \cite{irma-99}. A similar approach was later transformed to cactus graphs \cite{CaLSMaSi}, to outerplanar graphs \cite{MaSi12}, and to graphs with large enough girth \cite{CaFaSi,CaLiSi,koza-15}. For further reading about b-chromatic number and related concepts we recommend survey \cite{JaPe}.
        
There exist many variants of graph (vertex) colorings with some special extra condition(s). They usually yield a special chromatic number like acyclic chromatic number, star chromatic number, Thue chromatic number and many others. Their computational complexity is usually NP-hard and, similarly as in the case of chromatic number, we desire for some simple procedures that yield an upper bound for the mentioned invariants. Again, the information how much can go wrong in such a case is an interesting question. Therefore we start in this work with the analysis of acyclic b-chromatic number, that is the worst possible number of colors obtained by the second mentioned procedure which is limited in our case only to acyclic colorings of $G$.

The paper is organized as follows. In the next section we present basic notations and concepts, among others we recall two graph invariants: acyclic chromatic number and b-chromatic number. This part is followed by the definition of the acyclic b-chroatic number and some basic results about this parameter in Section \ref{sec_Definition}. Then we generalize an upper bound from the b-chromatic number to the acyclic b-chromatic number in Section \ref{secupperbound}. This allows us to present an upper bound that is quadratic with respect to the maximum degree of the graph. In Section \ref{sec_BAcyclicJoins} we present some results about the acyclic b-colorings of joins of graphs. The last section contains some final remarks and open problems.

%%%%%%%%%%%%%%%%%%%%%%%%%%%%%%%%%%%%%%%%%%%%%%%%%%%%%%%%%%%%%%%%%%%%%%%%%%%%%%%%

\section{Preliminaries}\label{sec_Preliminaries}

In this work we consider only graphs $G=(V(G),E(G))$ that are finite and simple, that is without loops and multiple edges. We use $n_G$ for the order and $m_G$ for the size of $G$. As usually we denote by $N_G(v)$ the \emph{open neighborhood} $\{u\in V(G):uv\in E(G)\}$ of $v\in V(G)$ and $N_G[v]=N_G(v)\cup\{v\}$ is the \emph{closed neighborhood} of $v$. The \emph{degree} of $v\in V(G)$ is denoted by $d_G(v)$ and is defined as $d_G(v)=\vert N_G(v)\vert$. By $\Delta (G)$ and $\delta (G)$ we denote the maximum and the minimum degree of a vertex in $G$, respectively. The clique number of $G$ is denoted by $\omega (G)$.  For $S\subseteq V(G)$ we denote by $G[S]$ the subgraph of $G$ induced by $S$. Graph $\overline{G}$ is the complement of $G$. We use $[k]$ to denote the set $\{1,\dots,k\}$ and $[j,k]$ to denote the set $\{j,\dots,k\}$ (so, in particular, $[1,k]=[k]$).

A graph $G$ is a \emph{cactus} graph if any two cycles intersect in at most one vertex. 
A graph $G$ is called an \emph{odd cycle graph} if $G$ does not contain any cycle of even length. By the best of our knowledge, we are not aware of the existence of the following basic result in the literature.

\begin{proposition}\label{oddcycle}
A graph $G$ is an odd cycle graph if and only if $G$ is a cactus graph with only odd cycles.
\end{proposition}

\noindent {\textbf{Proof.}} A cactus graph with only odd cycles is clearly an odd cycle graph by the definition. Otherwise, suppose that two cycles $C$ and $C'$ intersect in at least two vertices $u$ and $v$. We wish to show that there exists a cycle of even length. If $C$ or $C'$ is an even cycle, then we are done. So, we may assume that $C$ and $C'$ are odd cycles. We choose $u$ and $v$ in such a way that a $(u,v)$-path $P\subseteq C$ does not contain any vertex from $C'$ other than $u$ and $v$. Now $C^\prime$ splits into two $(u,v)$-paths $P_1$ and $P_2$ having no vertices in common with $P$ except their ends $u$ and $v$. Since $C'$ is odd, $P_1$ and $P_2$ have different parity. If $P_1$ has the same parity as $P$, then $P\cup P_1$ is an even cycle. Otherwise, $P_2$ has the same parity as $P$ and $P\cup P_2$ form an even cycle. So, $G$ is not an odd cycle graph.\qed\\

\subsection{Acyclic chromatic number}

A map $c:V(G)\rightarrow\{1,\ldots,k\}$ is called \emph{proper vertex coloring} with $k$ colors if $c(x)\neq c(y)$ for every edge $xy\in E(G)$. We consider here only proper vertex colorings, therefore we omit the terms "proper" and "vertex" and call $c$ a coloring or a $k$-coloring of $G$ in the remainder of the paper. The \emph{trivial coloring} of $G$ is the coloring where every vertex obtains a different color. The minimum number $k$, for which there exists a $k$-coloring, is called the \emph{chromatic number} of $G$ and is denoted by $\chi(G)$. Every $k$-coloring $c$ yields a partition of $V(G)$ into independent sets $V_i=\{u\in V(G):c(u)=i\}$, for every $i\in [k]$, called \emph{color classes} of $c$. We denote by $V_{i,j}$ the union $V_i\cup V_j$ and by $V_{i,j,\ell}$ the union $V_i\cup V_j\cup V_{\ell}$ for any $i,j,\ell\in [k]$. In particular we use $V_{i,j,\ell}(v)$ for component of $G[V_{i,j,\ell}]$ that contains vertex $v$. By $CN_c(v)$ we denote the set of all the colors that are present in $N_G(v)$ under coloring $c$, that is $CN_c(v)=\{c(u):u\in N_G(v)\}$. In addition we have $CN_c[v]=CN_c(v)\cup\{c(v)\}$.

A coloring $c$ is an \emph{acyclic coloring} of $G$ if $G[V_{i,j}]$ is a forest for any $i,j\in[k]$. In other words, the subgraph induced by any two color classes does not contain any cycle. Notice that $G[V_{i,i}]$ is even edgeless since $c$ is a coloring of $G$. The minimum number of colors of an acyclic coloring of $G$ is the \emph{acyclic chromatic number} denoted by $A(G)$. Clearly, $A(G)\geq \chi(G)$ as every acyclic coloring is also a coloring of $G$.

Acyclic colorings were introduced by B. Gr\"unbaum \cite{ref_Gru}, who proved that the acyclic chromatic number of any planar graph is not greater than $9$ and conjectured that in fact this bound is equal to $5$. This was finally proved by Borodin \cite{ref_Bor}. Mondal et al. \cite{ref_MonNisWhiRah} proved that every triangulated plane graph $G$ has an acyclically $3$-colorable subdivision, where the number of division vertices is not greater than $2.75n_G-6$, and an acyclically $4$-colorable subdivision, where the number of division vertices is not greater than $2n_G-6$. Alon, McDiarmid and Reed \cite{ref_AloMcDRee} showed that $A(G)\leq \lceil 50\Delta^{4/3}\rceil$, where $\Delta=\Delta(G)$, which proved the conjecture attributed to Erd\"os (see \cite[p.89]{ref_JenTof}), stating that $A(G)=o(\Delta^2)$. Recently, Gon\c{c}alves, Montassier and Pinlou \cite{ref_GonMonPin} used the entropy compression method to prove that for every graph $G$, $A(G)\leq \frac{3}{2}\Delta^{4/3}+O(\Delta)$. Alon, McDiarmid and Reed in \cite{ref_AloMcDRee} proved also that there exist graphs for which $A(G)=\Omega(\Delta^{4/3}/(\log\Delta)^{1/3})$. Alon, Mohar and Sanders \cite{ref_AloMohSan} showed that the acyclic chromatic number of the projective plane is $7$, while for every $G$ embeddable on a surface of Euler characteristic $\chi=-\gamma$, $A(G)=O(\gamma^{4/7})$. Moreover, for every $\gamma>0$ there exist graphs embeddable on surfaces of Euler characteristic $-\gamma$, for which $A(G)=\Omega(\gamma^{4/7}/(\log\gamma)^{1/7})$. This disproved the conjecture due to Borodin (see \cite[p.70]{ref_JenTof}) that acyclic chromatic number is equal to the chromatic number for all surfaces other than a plane. Acyclic colorings of graphs with bounded degree were studied e.g. by Fertin and Raspaud \cite{ref_FerRas}, Hocquard and Montassier \cite{ref_HocMon} and Yadav et al. \cite{ref_YadVarKotVen}. Mondal et al. \cite{ref_MonNisWhiRah} proved that deciding whether a graph with $\Delta\leq 6$ is acyclically $4$-colorable, is NP-complete.

\subsection{b-chromatic number}

Let $\mathcal{F}(G)$ be the set of all colorings of $G$ and let $c\in \mathcal{F}(G)$ be a $k$-coloring. A vertex $v$ of $G$ with $c(v)=i$ is a b-\emph{vertex} (of color $i$), if $v$ has all the colors of $c$ in its closed neighborhood, that is $CN_c[v]=[k]$. If a vertex $v$ with $c(v)=i$ is not a b-vertex, then (at least) one color is missing in $N_G[v]$, say $j$. We can recolor $v$ with $j$ and a slightly different coloring is obtained. Hence, if there exists no b-vertex of color $i$, then we can recolor every vertex $v$ colored with $i$ with some color not present in $N_G[v]$, say $j_v$. This way we obtain a new coloring $c_i:V(G)\rightarrow [k]\setminus\{i\}$ by 
\begin{equation*}
c_i(v)=\left\{ 
\begin{array}{ccc}
c(v) & : & c(v)\neq i \\
j_v & : & c(v)=i 
\end{array}
\right. .
\end{equation*}
Clearly $c_i$ is a $(k-1)$-coloring of $G$. We call the above procedure a  \emph{recoloring step}. By the \emph{recoloring algorithm} we mean an iterative performing of recoloring steps while it is possible where we start with a trivial coloring of $G$. Notice that one could also start with any other coloring from $\mathcal{F}(G)$.

Next we define the relation $\triangleleft$ on $\mathcal{F}(G)\times \mathcal{F}(G)$. We say that $c'\in \mathcal{F}(G)$ is in relation $\triangleleft$ with $c\in \mathcal{F}(G)$, $c'\triangleleft c$, if $c'$ can be obtained from $c$ by a recoloring step of one fixed color class of $c$.
Clearly, $\triangleleft$ is asymmetric. The transitive closure $\prec$ of $\triangleleft$ is then a strict partial order on $\mathcal{F}(G)$. Since there are finitely many different colorings of any graph $G$, this order has some minimal elements. The maximum number of colors of a minimal element of  $\prec$ is called the \emph{b-chromatic number} of $G$ and is denoted by $\varphi (G)$, in contrast with the chromatic number $\chi (G)$, which is the minimum number of colors of a minimal element of $\prec$. 

Every color class of every minimal element of $\prec$ needs to have a b-vertex. So, one can use the following alternative definition, as already mentioned in \cite{irma-99}. A coloring $c$ of $G$ is a \emph{b-coloring} if every color class contains a b-vertex. The b-chromatic number is then the maximum number of colors in a b-coloring of $G$. This definition was later used in almost all publications on the b-chromatic number. 

There exists a very natural upper bound for $\varphi(G)$. Namely, every b-coloring with $\varphi(G)$ colors needs at least $\varphi(G)$ vertices with $d_G(v)\geq\varphi(G)-1$. The $m$-\emph{degree} $m(G)$ is defined as 
$$m(G)=\max \{i: i-1\leq d_G(v_i)\},$$
where $v_1,\dots,v_{n_G}$ are ordered by the degrees $d_G(v_1)\geq \cdots\geq d_G(v_{n_G})$. It was shown already in \cite{irma-99} that $\varphi(G)\leq m(G)$. 

Every run of the recoloring algorithm gives a minimal element of $\prec$ and it is straightforward to see that it consists in at most $n_G-\chi(G)$ recoloring steps. In every recoloring step we need to find which color classes are without a b-vertex. If we use $\ell$ colors at some step of the recoloring algorithm, then only vertices of degree at least $\ell-1$ can be b-vertices (of a certain color). Hence we need to check only closed neighborhoods of such vertices and this can be done in $O(m_G)$ time in the worst case. Let us also mention that if a vertex $v$ with $c(v)=i$ is a b-vertex at some stage of the recoloring algorithm, then it remains a b-vertex of color $i$ after each recoloring step that still follows. Finally, when we have a color class $j$ without a b-vertex, we need to find for every vertex $v$ satisfying $c(v)=j$ a color that is not present in $N_G(v)$ and this can be done in $O(m_G)$ time again. Altogether, recoloring algorithm is polynomial algorithm and its time complexity is at most $O\left((n_G-\chi(G))(m_G)^2\right)$.    

With this we have a heuristic polynomial algorithm that produces a coloring of $G$ and gives an upper bound for $\chi(G)$. This way the study of $\varphi(G)$ became the study of the worst possible case that can happen while using the recoloring algorithm. A similar approach can be applied to the Grundy number $\Gamma(G)$, which corresponds with the study of the worst possible performance of the greedy algorithm. A comparative study on the b-chromatic number and the Grundy number was presented by Masih and Zaker \cite{MaZa1,MaZa2}. 
%%%%%%%%%%%%%%%%%%%%%%%%%%%%%%%%%
%%%%%%%%%%%%%%%%%%%%%%%%%%%%%%%%%

\section{Definition and some basic results}\label{sec_Definition}

Our goal in this part is to define the acyclic b-chromatic number of a graph $G$. For that reason we could be interested in colorings of $G$ that are acyclic and b-colorings at the same time. Unfortunately this is not possible for all graphs. Let us observe cycle $C_4$. Notice that $\varphi(C_4)=\chi(C_4)=2$ and the only b-coloring of $C_4$ is the 2-coloring that is not acyclic. So, both conditions, being an acyclic coloring and being a b-coloring,  are not always fulfilled. This is probably one of the reasons why this problem was not studied yet.  

We can avoid this problem if we focus more strictly on the original definition of $\varphi(G)$. For this let $\mathcal{AF}(G)$ be the set of all acyclic colorings of $G$. A recoloring step for $c\in\mathcal{AF}(G)$ that produces coloring $c'$ is an  \emph{acyclic recoloring step} if $c'\in \mathcal{AF}(G)$. This means that we can reduce some color of $c$ only when a new coloring $c'$ is also acyclic. The \emph{acyclic recoloring algorithm} is the use of an acyclic recoloring step until this is possible when starting by a trivial coloring of $G$. Also the acyclic recoloring algorithm has polynomial time complexity, because we need, in addition to the recoloring algorithm, after each acyclic recoloring step to check whether the new coloring is still acyclic. This can clearly be done in polynomial time since there are at most $O(n_G^2)$ pairs of different colors.

We define relation $\triangleleft_a\subseteq\mathcal{AF}(G)\times \mathcal{AF}(G)$ by $c'\triangleleft_a c$ when $c'$ can be obtained by an acyclic recoloring step from $c\in \mathcal{AF}(G)$. Similarly as $\triangleleft$, $\triangleleft_a$ is also asymmetric. Let $\prec_a$ be its transitive closure. Hence, $\prec_a$ is a strict partial order of $\mathcal{AF}(G)$. The trivial coloring is the greatest element of $\prec_a$ (sometimes also called the maximum element). Again, as $G$ is finite, also $\mathcal{AF}(G)$ is finite and at least one minimal element of $\prec_a$ exists.  

\begin{proposition}
 Let $G$ be a graph and $t\in\mathcal{AF}(G)$ be a trivial coloring. If $c\in\mathcal{AF}(G)$, then there exists a chain 
 $$c\triangleleft_a c_1\triangleleft_a c_2\triangleleft_a\cdots\triangleleft_a c_{\ell-1}\triangleleft_a t.$$
 \end{proposition}
 
\noindent {\textbf{Proof.}}
 Let $c$ be any $k$-coloring from $\mathcal{AF}(G)$. Let $\ell=n_G-k$, $c_{\ell}=t$ and $c_0=c$. We may assume that the $k$ colors from $c$ are the first $k$ colors from $t$. Let $v_1,\dots,v_k,v_{k+1},\dots,v_{n_G}$ be vertices of $G$ ordered in such a way that $c(v_i)=i$ for every $i\in [k]$, the rest of the ordering being arbitrary. For every $i\in [\ell]$ we define coloring $c_i$ from $c_{i-1}$ by 
 \begin{equation*}
c_i(v)=\left\{ 
\begin{array}{lcl}
c_{i-1}(v) & \text{if} & v\neq v_{k+i}, \\
k+i & \text{if} & v=v_{k+i}. 
\end{array}
\right.
\end{equation*}
In other words, if we reverse the order of colorings, then we obtain $c$ from $t$ by recoloring every vertex at most once. Clearly, $c_i\in\mathcal{AF}(G)$ for every $i\in [\ell]$ because $c\in\mathcal{AF}(G)$. Moreover, we can obtain $c_{i-1}$ from $c_i$ by a recoloring of vertex $v_{k+i}$. By the construction of $c_i$ from $c_{i-1}$ it is clear that $v_{k+i}$ is not a b-vertex for $c_i$ because the color $c_{i-1}(v_{k+i})$ is not in the closed neighborhood of $v_{k+i}$ in coloring $c_i$. Hence we have an acyclic recoloring step from $c_i$ to $c_{i-1}$  and $c_{i-1}\triangleleft_a c_i$ follows for every $i\in [\ell]$. \qed\bigskip

With this the following definition is justified. The \emph{acyclic b-chromatic number} $A_b(G)$ is the maximum number of colors in a minimal element of $\prec_a$:
$$A_b(G)=\max\{|c|:c\in\mathcal{AF}(G) \text{ is a minimal element of }\prec_a\}.$$
The acyclic b-chromatic number of a graph $G$ describes the worst case to appear while using the acyclic recoloring algorithm to estimate $A(G)$. An acyclic coloring of $G$ with $A_b(G)$ colors that arise from a minimal element of $\prec_a$ is called an $A_b(G)$-\emph{coloring}. We have the following inequality chain

\begin{equation}\label{chain}
\omega(G)\leq\chi(G)\leq A(G)\leq A_b(G)\leq n_G,
\end{equation}
where the fact that every acyclic b-coloring is also an acyclic coloring implies $A(G)\leq A_b(G)$. 
Next we characterize the graphs for which $A_b(G)=n_G$. 

\begin{proposition}\label{order}
We have $A_b(G)=n_G$ if and only if $G\cong K_{n_G}$. 
\end{proposition}

\noindent {\textbf{Proof.}}
If $G\cong K_n$, then $A_b(G)=n=n_G$ by (\ref{chain}). Conversely, if $G\ncong K_n$, then there exist different and nonadjacent $u,w\in V(G)$. For a trivial coloring $t$ of $G$,
 \begin{equation*}
c(v)=\left\{ 
\begin{array}{lccl}
t(v) &:& \text{if} & v\neq u, \\
t(w) &:& \text{if} & v=u, 
\end{array}
\right. 
\end{equation*}
is a coloring obtained by an acyclic recoloring step. Hence $t$ is not a minimal element of $\prec_a$ and $A_b(G)<n_G$.
\qed\bigskip

The b-chromatic number $\varphi(G)$ has an elegant description as a maximum number of colors for which b-vertex exists in every color class. For the acyclic b-chromatic number this is not enough as already shown for $C_4$. Namely, arbitrary recoloring can also trigger a bi-chromatic cycles and we need to avoid this. In order to formulate appropriate condition, we define the concept of a weak acyclic b-vertex.

\begin{definition}
Let $G$ be a graph with an acyclic coloring $c:V(G)\rightarrow [k]$. A vertex $v\in V_i$, $i\in [k]$, is a weak acyclic b-vertex if it satisfies
\begin{equation}\label{cond}
\forall \ell\in[k]-CN_c[v], \exists j\in CN_c(v): G[V_{j,\ell}\cup\{v\}] \text{ contains a cycle.}
\end{equation}
\end{definition}

Note that every b-vertex $v$ is also weak acyclic b-vertex, since $[k]-CN_c[v]=\emptyset$ in this case. If $v$ is a weak acyclic b-vertex and $\ell\in[k]-CN_c[v]$, then there exists a bi-colored path of even length between two neighbors of $v$ (both colored by $j$). For two colors $\ell,m\in[k]-CN_c[v]$ such two paths can have a common inner vertex of color $j$ if both paths start (and end) with the same color $j$.  

\begin{figure}[ht!]
\begin{center}
\begin{tikzpicture}[scale=0.5,style=thick,x=1cm,y=1cm]
\def\vr{3pt} % \vr = vertex radius;

% define vertices
%%%%%
%%%%%
\path (2,0) coordinate (a);
\path (0,0) coordinate (b);
\path (-1.5,1.5) coordinate (c);
\path (-1.5,3.5) coordinate (d);
\path (0,5) coordinate (e);
\path (2,5) coordinate (f);
\path (3.5,3.5) coordinate (g);
\path (3.5,1.5) coordinate (h);

%  edges

\draw (a) -- (b) -- (c) -- (d) -- (e) -- (f) -- (g) -- (h) -- (a);

\draw (a) [fill=white] circle (\vr);
\draw (b) [fill=white] circle (\vr);
\draw (c) [fill=white] circle (\vr);
\draw (d) [fill=white] circle (\vr);
\draw (e) [fill=white] circle (\vr);
\draw (f) [fill=white] circle (\vr);
\draw (g) [fill=white] circle (\vr);
\draw (h) [fill=white] circle (\vr);

\draw[anchor = north] (a) node {$1$};
\draw[anchor = north] (b) node {$3$};
\draw[anchor = east] (c) node {$1$};
\draw[anchor = east] (d) node {$3$};
\draw[anchor = south] (e) node {$1$};
\draw[anchor = south] (f) node {$2$};
\draw[anchor = west] (g) node {$1$};
\draw[anchor = west] (h) node {$2$};
\draw[anchor = south] (b) node {$y$};
\draw[anchor = west] (d) node {$x$};
\draw (1,2.5) node {$C_8$};

\path (10,0) coordinate (a1);
\path (8,0) coordinate (b1);
\path (6.5,1.5) coordinate (c1);
\path (6.5,3.5) coordinate (d1);
\path (8,5) coordinate (e1);
\path (10,5) coordinate (f1);
\path (11.5,3.5) coordinate (g1);
\path (11.5,1.5) coordinate (h1);

\path (6.5,0) coordinate (b2);
\path (6.5,5) coordinate (d2);
\path (13.5,3.5) coordinate (g2);
\path (13.5,1.5) coordinate (h2);
\path (13.5,0) coordinate (i);
%  edges

\draw (a1) -- (b1) -- (c1) -- (d1) -- (e1) -- (f1) -- (g1) -- (h1) -- (a1);
\draw (h1) -- (g2) -- (g1);
\draw (h1) -- (h2) -- (g1);
\draw (d1) -- (d2);
\draw (b1) -- (b2);
\draw (i) -- (h2);

\draw (a1) [fill=white] circle (\vr);
\draw (b1) [fill=white] circle (\vr);
\draw (c1) [fill=white] circle (\vr);
\draw (d1) [fill=white] circle (\vr);
\draw (e1) [fill=white] circle (\vr);
\draw (f1) [fill=white] circle (\vr);
\draw (g1) [fill=white] circle (\vr);
\draw (h1) [fill=white] circle (\vr);
\draw (b2) [fill=white] circle (\vr);
\draw (d2) [fill=white] circle (\vr);
\draw (g2) [fill=white] circle (\vr);
\draw (h2) [fill=white] circle (\vr);
\draw (i) [fill=white] circle (\vr);

\draw[anchor = north] (a1) node {$1$};
\draw[anchor = north] (b1) node {$3$};
\draw[anchor = east] (c1) node {$1$};
\draw[anchor = east] (d1) node {$3$};
\draw[anchor = south] (e1) node {$1$};
\draw[anchor = south] (f1) node {$2$};
\draw[anchor = east] (g1) node {$1$};
\draw[anchor = east] (h1) node {$2$};
\draw[anchor = south] (b1) node {$v$};
\draw[anchor = west] (d1) node {$u$};
\draw[anchor = west] (g2) node {$3$};
\draw[anchor = west] (h2) node {$4$};
\draw[anchor = east] (b2) node {$2(4)$};
\draw[anchor = north] (b2) node {$z$};
\draw[anchor = east] (d2) node {$4$};
\draw[anchor = south] (g2) node {$w$};
\draw[anchor = west] (i) node {$3$};
\draw[anchor = east] (i) node {$t$};
\draw[anchor = north] (h1) node {$b$};
\draw[anchor = south] (g1) node {$a$};
\draw (13,1) node{$c$};
\draw (9,2.5) node {$G$};
\end{tikzpicture}
\end{center}
\caption{Colorings without weak acyclic b-vertices of all colors that can sometimes be acyclic recolored and sometimes not.}
\label{notenough}
\end{figure}
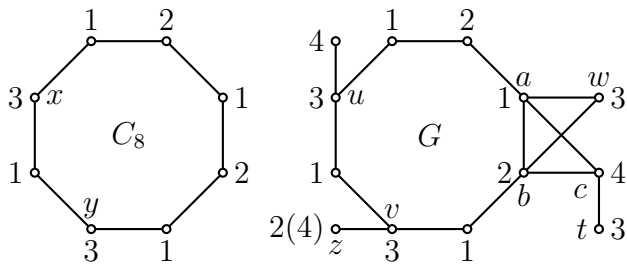

As already suggested by "weak", the notion of weak acyclic b-vertex does not generalize b-vertices to acyclic b-vertices in all the cases. For this observe Figure \ref{notenough}. On the left we have a colored $8$-cycle where only color 1 has a weak acyclic b-vertex (that is actually a b-vertex), while colors $2$ and $3$ have no weak acyclic b-vertex. However, the presented coloring cannot be reduced with an acyclic recoloring step. Even more tricky is the example on the right graph $G$ of Figure \ref{notenough}. There are two similar colorings of this graph with only difference in vertex $z$. In first coloring $z$ is colored by $2$ and in the second coloring by $4$. Vertices $a,b,c$ are b-vertices for colors $1,2,4$, respectively, so only color $3$, which is without weak acyclic b-vertex, can be recolored. In first coloring, this is possible as $w$ and $v$ can be recolored with $4$, and $u$ and $t$ with $2$ and we get an acyclic coloring. On the other hand, an acyclic recoloring of $G$ is not possible for the second coloring (i.e. the one in which $c(z)=4$), because $u$ and $v$ can be recolored only with $2$ and this yields a bi-chromatic cycle. 

\begin{figure}[ht!]
\begin{center}
\begin{tikzpicture}[scale=0.5,style=thick,x=1cm,y=1cm]
\def\vr{3pt} % \vr = vertex radius;

% define vertices
%%%%%
%%%%%
\path (2,0) coordinate (a);
\path (0,0) coordinate (b);
\path (-1.5,2.5) coordinate (c);
\path (0,5) coordinate (e);
\path (2,5) coordinate (f);
\path (3.5,2.5) coordinate (g);
\path (5,0) coordinate (d);
\path (7,0) coordinate (h);
\path (9,0) coordinate (i);
\path (10.5,2.5) coordinate (j);
\path (5,5) coordinate (k);
\path (7,5) coordinate (l);
\path (9,5) coordinate (m);

\path (-1,7) coordinate (e1);
\path (1,7) coordinate (e2);
\path (12,1.5) coordinate (j1);
\path (12,3.5) coordinate (j2);
\path (2,-2) coordinate (a1);
\path (5,-2) coordinate (a2);
\path (0,-2) coordinate (b1);
%  edges

\draw (g) -- (a) -- (b) -- (c) -- (e) -- (e) -- (f) -- (g) -- (d) -- (b) -- (b1);
\draw (b) -- (a1) -- (a2) -- (a) -- (a1) -- (d) -- (a2);
\draw (a) -- (d) -- (h) -- (i) -- (j) -- (m) -- (l) -- (k) -- (g) -- (d);
\draw (e1) -- (e) -- (e2);
\draw (j1) -- (j) -- (j2);

\draw (a) [fill=white] circle (\vr);
\draw (b) [fill=white] circle (\vr);
\draw (c) [fill=white] circle (\vr);
\draw (d) [fill=white] circle (\vr);
\draw (e) [fill=white] circle (\vr);
\draw (f) [fill=white] circle (\vr);
\draw (g) [fill=white] circle (\vr);
\draw (h) [fill=white] circle (\vr);
\draw (i) [fill=white] circle (\vr);
\draw (j) [fill=white] circle (\vr);
\draw (k) [fill=white] circle (\vr);
\draw (l) [fill=white] circle (\vr);
\draw (m) [fill=white] circle (\vr);
\draw (a1) [fill=white] circle (\vr);
\draw (a2) [fill=white] circle (\vr);
\draw (e1) [fill=white] circle (\vr);
\draw (j1) [fill=white] circle (\vr);
\draw (j2) [fill=white] circle (\vr);
\draw (e2) [fill=white] circle (\vr);
\draw (b1) [fill=white] circle (\vr);

\draw[anchor = south] (a) node {$2$};
\draw[anchor = south] (b) node {$1$};
\draw[anchor = east] (c) node {$2$};
\draw[anchor = south] (d) node {$4$};
\draw[anchor = north] (e) node {$3$};
\draw[anchor = south] (f) node {$2$};
\draw[anchor = west] (g) node {$1$};
\draw[anchor = north] (h) node {$5$};
\draw[anchor = north] (b1) node {$4$};
\draw[anchor = north] (a1) node {$5$};
\draw[anchor = north] (a2) node {$3$};
\draw[anchor = south] (e1) node {$4$};
\draw[anchor = south] (e2) node {$5$};
\draw (1,2.5) node {$C$};

\draw[anchor = north] (i) node {$4$};
\draw[anchor = north] (j) node {$1$};
\draw[anchor = south] (k) node {$4$};
\draw[anchor = south] (l) node {$5$};
\draw[anchor = south] (m) node {$4$};
\draw[anchor = south] (j1) node {$2$};
\draw[anchor = south] (j2) node {$3$};
\draw (7,2.5) node {$C'$};

\draw[anchor = east] (b) node {$a$};
\draw (1.5,-0.5) node {$b$};
\draw (5.5,-0.5) node {$c$};
\draw[anchor = east] (g) node {$d$};
\draw[anchor = east] (a1) node {$g$};
\draw[anchor = east] (e) node {$e$};
\draw[anchor = east] (b1) node {$h$};
\draw[anchor = south] (j) node {$f$};
\end{tikzpicture}
\end{center}
\caption{Cycles $C$ and $C'$ form a critical cycle system.}
\label{system}
\end{figure}
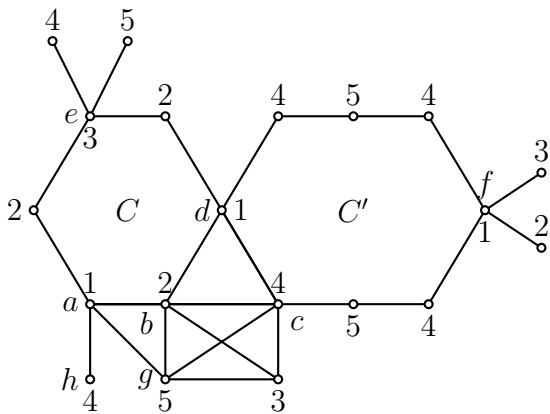

Both colorings of $G$ in Figure \ref{notenough} show that sometimes also an even cycle $C$ can prevent the acyclic recoloring step for a color $i$ and not a weak acyclic b-vertex. This is possible only if $C$ is colored with exactly three colors, say $i,j,k$, and $i$ appears at least twice on $C$. Further, one of the other two colors, say $j$, must be on every second vertex. We call such a cycle an $i$-\textit{critical cycle} or $i$-CC for short. Note, that color $k$ can appear only once, say on vertex $v$, on $C$ and in such a case $v$ cannot be recolored with $i$ by condition (\ref{cond}) (when $v$ is a weak b-acyclic vertex). However, if also $k$ appears at least twice on $C$, then $C$ is $k$-CC as well and we say that $C$ is $i,k$-critical cycle or $i,k$-CC for short. Color $i$ of $i$-CC is called \textit{principal} while in $i,k$-CC we have two principal colors $i$ and $k$. Clearly, $i$-CC is of the length at least six and $i,k$-CC of the length at least eight. Note that $C_8$ from Figure \ref{notenough} is $2,3$-CC and the $8$-cycle from $G$ on the same figure is $2,3$-CC as well. Further, in Figure \ref{system} $C$ is $1$-CC and $C'$ is $1,5$-CC.

Critical cycles that have a common vertex of principal color can further influence on each other. For this observe a simple example in Figure \ref{system}. Here $C$ is $1$-CC and $C'$ is $1,5$-CC and they have a common vertex $d$ of principal color $1$. Notice that $b,c,g$ are b-vertices of colors $2,4,5$, respectively. In addition, $e$ is a weak acyclic b-vertex. So the only candidate for acyclic recoloring is color $1$. In this case notice that $f$ can be recolored only with $5$. Now, $d$ can be recolored with $3$ or $5$, but we cannot recolor both $d$ and $f$ with $5$, because then $C'$ is bi-colored. So, $d$ must receive color $3$. Finally, $a$ can be recolored only with $3$, which yields bi-colored cycle $C$. Hence, acyclic recoloring of this coloring is not possible.

Let $c$ be an acyclic coloring of $G$. Let $\mathcal{C}$ be a collection of all critical cycles for $c$. A \textit{critical cycle system for color $i$} or CCS$(i)$ for short is a subcollection of $i$-CC from $\mathcal{C}$ that form a maximal connected subgraph of $G$ and where different critical cycles intersect only in vertices of color $i$ (there can be more than one CCS$(i)$ in one coloring of $G$). The color that is not principal on a critical cycle $C$ appears exactly once on $C$ or on half of vertices of $C$. If it appears only once, then the condition (\ref{cond}) prevents the recoloring with $i$ if necessary. The other color, say $j$, can be always recolored on $C$ (if there is no weak acyclic b-vertex of this color or some other CCS$(j)$). In Figure \ref{system} $\{C,C'\}$ is the collection of all $1$-CC as well as the only CCS$(1)$ for the presented coloring. On the other hand $C'$ is the only $5$-CC and with this also CCS$(5)$.

Let $i$ be a principal color without any weak acyclic b-vertex in a coloring of a graph $G$. Let $v$ be a vertex of color $i$ in a CCS$(i)$. An \textit{available color} for $v$ is every color that is not in the neighborhood of $v$ and is not on a bi-colored path of even length between two neighbors of $v$. By $A_v$ we denote the set of all available colors of $v$. In Figure \ref{system} we have $A_a=\{3\}$, $A_d=\{3, 5\}$ and $A_f=\{5\}$ for CCS$(1)$. On the same figure for CCS$(5)$, both sets of available colors are empty, because there exists a weak acyclic b-vertex $g$ of color $5$.  

Let $D\subseteq \mathcal{C}$ be a CCS$(i)$ for an acyclic coloring of graph $G$. If there exists a principal color $i$ among the colors used on $D$ such that for every vertex $v$ from $D$ colored with $i$ there exists a color $j_v\in A_v$, such that recoloring all such $v$ with $j_v$ results in an acyclic coloring, then $D$ is \textit{recolorable}. Otherwise, $D$ is not recolorable. CCS$(1)$ $\{C,C'\}$ from Figure \ref{system} is not recolorable. As we have already observed, the recoloring of color $1$ leads to a bi-chromatic cycle. Vertices of color $5$ could be recolored in CCS$(5)$, however vertex $g$ is a b-vertex of color $5$ and with this also a weak acyclic b-vertex and color $5$ cannot be recolored. It is also easy to see that the critical cycle $C_8$ from Figure \ref{notenough}  is not recolorable as well as the second coloring of $G$ on the same figure. However, the first coloring of this graph is recolorabe. Now everything is settled for the definition of acyclic b-vertices.

\begin{definition}
Let $G$ be a graph with an acyclic coloring $c:V(G)\rightarrow [k]$. A vertex $v\in V_i$, $i\in [k]$, is an acyclic b-vertex if it satisfies
\begin{equation}\label{cond1}
\forall \ell\in[k]-CN_c[v], \exists j\in CN_c(v): (G[V_{j,\ell}\cup\{v\}] \text{ contains a cycle }\vee
\end{equation}
\begin{center}
there exists a CCS$(i)$ of $G$ that contains $v$ and is not recolorable).
\end{center}
\end{definition}

Now we can describe minimal elements of $\prec_a$ of graph $G$ as follows.

\begin{theorem}\label{conditionABcoloring}
An acyclic $k$-coloring $c$ is a minimal element of $\prec_a$ if and only if every color class $V_i$, $i\in[k]$, contains an acyclic b-vertex.
\end{theorem}

\noindent {\textbf{Proof.}}
Let an acyclic $k$-coloring $c$ be a minimal element of $\prec_a$ of a graph $G$. This means that we cannot present an acyclic recoloring step for $c$. There are two possible reasons for that for any color class $V_i$, $i\in [k]$. Firstly, $V_i$ has a b-vertex or secondly, by any recoloring of $V_i$ we get a bi-chromatic cycle of colors, say $j$ and $\ell$, different than $i$. Let $v_1v_2\dots v_qv_1$ be such a cycle $C$ where $v_1\in V_i$ and $c(v_2)=c(v_q)=j$. Clearly, $C$ can be different after different recolorings. Since the recoloring is arbitrary, we may assume that $\ell\notin CN_c[v_1]$ is arbitrary. If $v_1$ is the only vertex of color $i$ on $C$, then path $v_2\dots v_q$ is bi-colored by $c$. Hence, $G[V_{j,\ell}\cup\{v\}$ contains a cycle and the first part of condition (\ref{cond1}) holds. Otherwise, there exists more vertices of color $i$ on $C$. Then there exists CCS$(i)$ that contains $v_1$, because $C$ contains more than one vertex of color $i$ and some $i$-CC from the mentioned CCS$(i)$ is bi-chromatic after any recoloring. So, there exists CCS$(i)$ that is not recolorable and the second part of condition (\ref{cond1}) is fulfilled. In both cases we have an acyclic b-vertex for every color $i\in [k]$. 

Conversely, if $c$ is not a minimal element of $\prec_a$, then we can perform an acyclic  recoloring step for some color $i$. This means that for every $v\in V_i$ there exists a color $\ell_v\notin CN_c[v]$ such that for every $j\in CN_c(v)$ there is no cycle in $G[V_{j,\ell_v}\cup\{v\}]$ and every CCS$(i)$ (if it exists) is acyclic recolorable. But then condition (\ref{cond1}) is not fulfilled and $V_i$ is without acyclic b-vertex and we are done. 
\qed\bigskip

\begin{corollary}\label{cor}
The acyclic b-chromatic number $A_b(G)$ of a graph $G$ is the largest integer $k$, such that there exists an acyclic $k$-coloring, where every color class $V_i$, $i\in [k]$, contains an acyclic b-vertex.
\end{corollary}

\begin{corollary}\label{cor1}
Let $G$ be a graph with all even cycles being $4$-cycles. The acyclic b-chromatic number $A_b(G)$ of $G$ is the largest integer $k$, such that there exists an acyclic $k$-coloring, where every color class $V_i$, $i\in [k]$, contains a weak acyclic b-vertex.
\end{corollary}

\noindent{}Let us observe some simple facts.

\begin{corollary}\label{basic2}
For every positive integers $n,k,\ell$, where $k\geq 3$ and $\ell\geq 5$,  we have
\begin{itemize}
\item $A_b(\overline{K}_n)=1$.
\item $A_b(P_{\ell})=3$.
\item $A_b(C_k)=3$.
\end{itemize}
\end{corollary}

\noindent {\textbf{Proof.}}
$A_b(\overline{K}_n)\geq 1$ and $A_b(C_k)\geq 3$ follow from (\ref{chain}), while the coloring $c:V(P_{\ell})\rightarrow [3]$ guaranteeing $A_b(P_{\ell})\geq 3$ can be defined as $c(v_i)=(i\imod 3$)+1, $i\in [\ell]$ (every internal vertex is a b-vertex and thus there is an acyclic b-vertex in each color class). On the other hand, if $\overline{K}_n$ is colored with $p\geq 2$ colors, then every color is without an acyclic b-vertex and this coloring is not a minimal element of $\prec_a$ by Theorem \ref{conditionABcoloring}. Similar is with $C_k$ and $P_{\ell}$. If they are colored by $p\geq 4$ colors, then no color class contains an acyclic b-vertex and such a coloring in not a minimal element of $\prec_a$ by Theorem \ref{conditionABcoloring}. The desired equalities now follow.
\qed\bigskip

We end this section with a brief discussion on graphs where $\mathcal{F}(G)=\mathcal{AF}(G)$, or in other words, where every coloring of $G$ is acyclic. Among such graphs are clearly odd cycle graphs described in Proposition \ref{oddcycle}, some cactus graphs and in particular trees. For such graphs we have $A_b(G)=\varphi(G)$. In particular, b-chromatic number of trees and cactus graphs was studied in \cite{irma-99} and \cite{CaLSMaSi}, respectively. In both cases it was shown that 
$$m(G)-1\leq\varphi(G)\leq m(G),$$
where the above holds for cactus graphs when $m(G)\geq 7$. 
Moreover, the lower bound is achieved if and only if $G$ is a pivoted graph. Notice that pivoted tree (see \cite{irma-99}) is defined differently than a pivoted cactus graph (see \cite{CaLSMaSi}) and that for pivoted trees we do not have a restriction that $m(G)\geq 7$. It is also important to mention, that the authors of \cite{CaLSMaSi} used an alternative definition of cacti, where two cycles can have arbitrarily many vertices in common, provided that they do not have a common edge. For that reason their results are not consistent with the definitions used in our paper. However, we still can formulate the following.

\begin{corollary}\label{tree}
Let $T$ be a tree. If $T$ is a pivoted tree, then $A_b(T)=m(T)-1$ and otherwise $A_b(T)=m(T)$.
\end{corollary}

Corollary \ref{tree} implies in particular that the difference between $A_b(G)$ and $A(G)$ can be arbitrarily large.
\begin{corollary}\label{treeInfinite}
There exists an infinite family of graphs $G_1, G_2, \dots$ such that $(A_b(G_n)-A(G_n))\rightarrow\infty$ as $n\rightarrow\infty$. 
\end{corollary}
\noindent {\textbf{Proof.}}
For $n\geq 1$, let $G_n$ be a graph consisting of a star $K_{1,n+1}$ with $n$ pendant edges attached to each of its leaves. Graph $G_n$ has exactly $n+2$ vertices of degree $n+1$ and $n(n+1)$ vertices of degree $1$, so $m(G_n)=n+2$. Also, $G_n$ is not a pivoted tree, so $A_b(G_n)=m(G_n)=n+2$. On the other hand, any proper coloring of $G_n$ is its acyclic coloring, so $A(G)=2$. Thus $A_b(G_n)-A(G_n)=n$ for every $n\geq 1$ and $(A_b(G_n)-A(G_n))\rightarrow\infty$ as $n\rightarrow\infty$. 
\qed\bigskip

%%%%%%%%%%%%%%%%%%%%%%%%%%%%%%%%%%%%%%%%%%%%%%%%%%%%%%%%%%
%%%%%%%%%%%%%%%%%%%%%%%%%%%%%%%%%%%%%%%%%%%%%%%%%%%%%%%5

\section{An upper bound on $A_b(G)$ analogous to $m(G)$ for $\varphi(G)$}\label{secupperbound}

Let $c$ be an acyclic $k$-coloring of a graph $G$ that is minimal with respect to order $\prec_a$. Recall that according to Theorem \ref{conditionABcoloring} every color class $V_i$ contains an acyclic b-vertex. While for b-vertices, high enough degree is necessary, for acyclic b-vertices we need high enough degree or sufficient number of even-vertex internally (or EVI for short) disjoint $(u,w)$-paths of odd order (that is, of even length) between some of its neighbors $u\neq w$ or a combination of both (by even-vertex internally disjoint we mean that two such paths can have odd vertices of such a path in common, but not the even vertices). In particular, if for $u,w\in N_G(v)$ there are $k$ different EVI disjoint $(u,w)$-paths $P_1,\dots,P_k$ of even length, where $P_k=uvw$, in the worst case these paths can yield $k$ different colors into which $v$ cannot be recolored. Indeed, on every path $P_i$, $i\in[k-1]$, there can exist $a_i\in V(P_i)$ with $c(a_i)\notin\{c(v),c(u)\}$ and we can have alternating colors $c(u)=c(w)$ and $c(a_i)$ or $P_k\cup P_i$ is $c(v)$-CC that contains colors $c(v),c(u)$ and $c(a_i)$ and belongs to a CCS$(c(v))$ that  is not recolorable. Hence, in the worst case, $u, a_1,\dots ,a_{k-1}$ are colored differently and $v$ cannot receive any of their colors in an acyclic recoloring step.

This can be generalized to a bigger number of neighbors in the following way. Consider a weak partition $P=\{A^P_0,A^P_1,\dots,A^P_k\}$ of $N_G(v)$ into $k+1$ disjoint sets such that $|A^P_0|\geq 0$ and $|A^P_i|\geq 2$ for $i\in [k]$ (the mentioned partition is weak because $A^P_0$ can be empty). Without loss of generality assume that $v$ is colored with color $1$, the vertices of $A^P_0$ with distinct colors from the set $[2,|A^P_0|+1]$ and all the vertices of $A^P_i$, $i\in [k]$, with color $|A^P_0|+i+1$. Now, let ${\rm elp}_G(v,P)$ be the maximum number of pairwise EVI disjoint paths disjoint with $v$ having odd number of vertices (i.e., of even length), with both ends in one of the sets $A^P_i$, $i\in [k]$. Observe that given a partition with color classes defined as above, in the worst case one cannot recolor $v$ to exactly $(|A^P_0|+k+{\rm elp}_G(v,P))$ colors different than $c(v)$: $(|A^P_0|+k)$ colors are blocked by the neighbors and ${\rm elp}_G(v,P)$ by the alternately colored bi-chromatic EVI disjoint paths that could appear in the coloring. This encourages us to define the \textit{acyclic degree} of $v$ as
$$
d_G^a(v)=\max_{P\in\mathcal{P}(v)}\{(|A^P_0|+(|P|-1)+{\rm elp}_G(v,P))\},
$$  
where $\mathcal{P}(v)$ is the family of all the weak partitions $P$ of $N_G(v)$ defined as above.

\begin{figure}[ht!]
\begin{center}
\begin{tikzpicture}[scale=1,style=thick,x=1cm,y=1cm]
\def\vr{2pt} % \vr = vertex radius;

% define vertices
%%%%%
%%%%%
\path (2,0) coordinate (a);
\path (2,4) coordinate (b);
\path (1.5,2) coordinate (c);
\path (0,2) coordinate (d);
\path (2.5,2) coordinate (e);
\path (4,2) coordinate (f);
\path (5.5,2) coordinate (g);
\path (5.3,3) coordinate (h);
\path (5.3,1) coordinate (i);

%  edges

\draw (g) -- (f) -- (b) -- (e) -- (a) -- (d) -- (b) -- (c) -- (a) -- (f);
\draw (h) -- (f) -- (i);

\draw (a) [fill=white] circle (\vr);
\draw (b) [fill=white] circle (\vr);
\draw (c) [fill=white] circle (\vr);
\draw (d) [fill=white] circle (\vr);
\draw (e) [fill=white] circle (\vr);
\draw (f) [fill=white] circle (\vr);
\draw (g) [fill=white] circle (\vr);
\draw (h) [fill=white] circle (\vr);
\draw (i) [fill=white] circle (\vr);

\draw[anchor = north] (a) node {$x_1^1$};
\draw[anchor = south] (b) node {$x_2^1$};
\draw[anchor = east] (c) node {$y_3^1$};
\draw[anchor = east] (d) node {$y_4^1$};
\draw[anchor = west] (e) node {$y_2^1$};
\draw[anchor = south] (f) node {$y_1^1$};
\draw[anchor = west] (g) node {$z_1^1$};
\draw[anchor = west] (h) node {$y_1^2$};
\draw[anchor = west] (i) node {$u$};

\end{tikzpicture}
\end{center}
\caption{Graph $G$ with the optimal weak partition $A^P_0=\{u,z_1^1,y_1^2\}$,$A^P_1=\{x_1^1,x_2^1\}$, implying $|A^P_0|=3$, $|P|-1=1$, ${\rm elp}_G(v,P)=3$ and $d_G^a(y_1^1)=7$.}
\label{elp}
\end{figure}
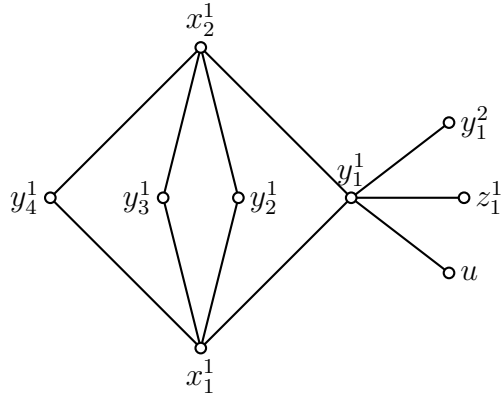

See Figure \ref{elp} with the only optimal weak partition of $N_G(y_1^1)$ defined as $P=\{A^P_0,A^P_1\},A^P_0=\{u,z_1^1,y_1^2\}$,$A^P_1=\{x_1^1,x_2^1\}$. Clearly, $|A^P_0|=3$, $|P|-1=1$ and ${\rm elp}_G(v,P)=3$ (the paths being $x_1^1y_2^1x_2^1, x_1^1y_3^1x_2^1, x_1^1y_4^1x_2^1$), thus $d_G^a(y_1^1)=7$.

Since we traverse over all weak partitions in $\mathcal{P}(v)$, notice that $d_G^a(v)$ represents the maximum number of colors in $N_G(v)$ and on  EVI disjoint paths of even length between vertices of $N_G(v)$ into which $v$ cannot be recolored in a recoloring step. This gives an analogy to the relation between degree $d_G(v)$ and $\varphi(G)$, where we needed sufficient number of vertices of high degree to expect $\varphi(G)$ to be large. Hence, we are encouraged to define an $m_a$-degree of a graph $G$ denoted by $m_a(G)$. First we order the vertices $v_1,\dots,v_{n_G}$ of $G$ by non-increasing acyclic degree. The value of $m_a(G)$ is then the maximum position $i$ in this order such that $d_G^a(v_i)\geq i-1$, that is
$$m_a(G)=\max \{i: i-1\leq d_G^a(v_i)\}.$$

\begin{theorem}\label{madegree}
For any graph $G$ we have $A_b(G)\leq m_a(G)$.
\end{theorem}

\noindent {\textbf{Proof.}}
On the way to a contradiction suppose that there exists a graph $G$ for which $k=A_b(G)>m_a(G)$ and that $c:V(G)\rightarrow [k]$ is an  appropriate acyclic b-coloring of $G$. We may assume that $v_1,\dots, v_{n_G}$ are ordered by non-increasing acyclic degree. Clearly, not all colors from $[k]$ are present on the vertices $v_1,\dots,v_{k-1}$ and we may assume that $c(v_i)\neq k$ for every $i\in[k-1]$. This means that $d_G^a(v)<k-1$ for every vertex $v$ of color $k$. Hence, for every vertex $v$ of color $k$, (\ref{cond1}) does not hold, a contradiction with Theorem \ref{conditionABcoloring}. 
\qed\bigskip

In general case it seems to be hard to compute $m_a(G)$, because we need to derive acyclic degrees for every vertex $v$ of $G$. This is further connected with all the weak partitions from $\mathcal{P}(v)$. Moreover, for every such weak partition we need to get the maximum number of EVI disjoint paths of even length. 

The problem of finding a maximum cardinality set of disjoint paths between two fixed vertices is obviously related to the connectivity of graph and Menger's Theorem \cite{ref_Menger} and has been studied for many years. It was proved to be polynomially solvable already by Even and Tarjan \cite{ref_EvenTarjan} with network flow algorithms and no spectacular progress has been made since then, although some papers may be found giving better results in special cases (see e.g. \cite{ref_GrossiMarinoVersari} or \cite{ref_PreisserSchmidt} for some recent results). Unfortunately, when the path lengths are somehow restricted, the problem was proved to be NP-complete already by Itai et al. \cite{ref_ItaiPerlShiloach}. Bley \cite{ref_Bley} proved its APX-completeness. Also the version of the problem of finding $k$ disjoint paths between $k$ disjoint pairs of terminals has been studied, see e.g. Fleszar et al. \cite{ref_FleszarMnichSpoerhase} for some recent results. We do not know any results about the problem of finding maximum cardinality set of EVI disjoint paths between arbitrary vertices, even if they are members of one fixed set.

Nevertheless, one can expect some success with $m_a(G)$ for certain graph classes, mainly the ones with low connectivity which are not too dense or with low maximum degree. In the next part of this section we demonstrate this approach on a family of graphs. It will allow us to show that there is no linear relationship between $A_b(G)$ and $\Delta(G)$ in general case. This underlines another difference between $A_b(G)$ and $\varphi (G)$. Recall that $\varphi(G)\leq m(G)\leq \Delta(G)+1$, however, as we are going to show, $A_b(G)$ can be arbitrarily larger than $\Delta(G)$.

\begin{theorem}\label{roofsInfinite}
There exists an infinite family of graphs $G_1, G_2, \dots$ such that $(A_b(G_n)-\Delta(G_n))\rightarrow\infty$ as $n\rightarrow\infty$. 
\end{theorem}
\noindent {\textbf{Proof.}}
For $n\geq 1$, let $H_n^i$, $i\in[2n+4]$, be isomorphic graphs with 
$$V(H_n^i)=\{y_1^i,\dots,y_{n+2}^i\}\cup\{x_1^i,x_2^i\}\cup\{z_1^i,\dots,z_{n-1}^i\}$$
and 
$$E(H_n^i)=\{x_j^iy_{\ell}^i:j\in[2],\ell\in[n+2]\}\cup\{y_1^iz_j^i:j\in[n-1]\}.$$
Notice that in Figure \ref{elp} a graph $H_2^1$ is presented with two additional leaves $u$ and $y_1^2$. Using graphs $H_n^i$, $i\in[2n+4]$, we construct graph $G_n$ with 
$$V(G_n)=\{u,v\}\cup \bigcup_{i=1}^{2n+4}V(H_n^i)$$
and
$$E(G_n)=\{uy_1^1,vy_1^{2n+4},y_1^iy_1^{i+1}:i\in[2n+3]\}\cup \bigcup_{i=1}^{2n+4}E(H_n^i).$$
Observe $G_2$ in Figure \ref{G_5}, while the first part $H_2^1$ together with the neighbors $u$ and $y_1^2$ of $y_1^1$ in $G_2$ is presented in Figure \ref{elp}.

\begin{figure}[ht!]
\begin{center}
\begin{tikzpicture}[scale=1.5,style=thick,x=1cm,y=1cm]
\def\vr{2pt} % \vr = vertex radius;

% define vertices
%%%%%
%%%%%
\path (0,0) coordinate (q0);
\path (1,0) coordinate (q1);
\path (2,0) coordinate (q2);
\path (3,0) coordinate (q3);
\path (4,0) coordinate (q4);
\path (5,0) coordinate (q5);
\path (6,0) coordinate (q6);
\path (7,0) coordinate (q7);
\path (8,0) coordinate (q8);
\path (9,0) coordinate (q9);

\path (1,-1) coordinate (p1);
\path (2,-1) coordinate (p2);
\path (3,-1) coordinate (p3);
\path (4,-1) coordinate (p4);
\path (5,-1) coordinate (p5);
\path (6,-1) coordinate (p6);
\path (7,-1) coordinate (p7);
\path (8,-1) coordinate (p8);

\path (0.65,1) coordinate (r1);
\path (1.65,1) coordinate (r2);
\path (2.65,1) coordinate (r3);
\path (3.65,1) coordinate (r4);
\path (4.65,1) coordinate (r5);
\path (5.65,1) coordinate (r6);
\path (6.65,1) coordinate (r7);
\path (7.65,1) coordinate (r8);

\path (1.35,1) coordinate (s1);
\path (2.35,1) coordinate (s2);
\path (3.35,1) coordinate (s3);
\path (4.35,1) coordinate (s4);
\path (5.35,1) coordinate (s5);
\path (6.35,1) coordinate (s6);
\path (7.35,1) coordinate (s7);
\path (8.35,1) coordinate (s8);

\path (1,1.5) coordinate (t1);
\path (2,1.5) coordinate (t2);
\path (3,1.5) coordinate (t3);
\path (4,1.5) coordinate (t4);
\path (5,1.5) coordinate (t5);
\path (6,1.5) coordinate (t6);
\path (7,1.5) coordinate (t7);
\path (8,1.5) coordinate (t8);

\path (1,2.5) coordinate (u1);
\path (2,2.5) coordinate (u2);
\path (3,2.5) coordinate (u3);
\path (4,2.5) coordinate (u4);
\path (5,2.5) coordinate (u5);
\path (6,2.5) coordinate (u6);
\path (7,2.5) coordinate (u7);
\path (8,2.5) coordinate (u8);

\path (1,3.5) coordinate (v1);
\path (2,3.5) coordinate (v2);
\path (3,3.5) coordinate (v3);
\path (4,3.5) coordinate (v4);
\path (5,3.5) coordinate (v5);
\path (6,3.5) coordinate (v6);
\path (7,3.5) coordinate (v7);
\path (8,3.5) coordinate (v8);

%  edges
\draw (q0) -- (q9);
\draw (q1) -- (p1);
\draw (q2) -- (p2);
\draw (q3) -- (p3);
\draw (q4) -- (p4);
\draw (q5) -- (p5);
\draw (q6) -- (p6);
\draw (q7) -- (p7);
\draw (q8) -- (p8);

\draw (q1) -- (r1) -- (u1) -- (s1);
\draw (q2) -- (r2) -- (u2) -- (s2);
\draw (q3) -- (r3) -- (u3) -- (s3);
\draw (q4) -- (r4) -- (u4) -- (s4);
\draw (q5) -- (r5) -- (u5) -- (s5);
\draw (q6) -- (r6) -- (u6) -- (s6);
\draw (q7) -- (r7) -- (u7) -- (s7);
\draw (q8) -- (r8) -- (u8) -- (s8);

\draw (r1) -- (v1) -- (s1);
\draw (r2) -- (v2) -- (s2);
\draw (r3) -- (v3) -- (s3);
\draw (r4) -- (v4) -- (s4);
\draw (r5) -- (v5) -- (s5);
\draw (r6) -- (v6) -- (s6);
\draw (r7) -- (v7) -- (s7);
\draw (r8) -- (v8) -- (s8);

\draw (q1) -- (s1) -- (t1) -- (r1);
\draw (q2) -- (s2) -- (t2) -- (r2);
\draw (q3) -- (s3) -- (t3) -- (r3);
\draw (q4) -- (s4) -- (t4) -- (r4);
\draw (q5) -- (s5) -- (t5) -- (r5);
\draw (q6) -- (s6) -- (t6) -- (r6);
\draw (q7) -- (s7) -- (t7) -- (r7);
\draw (q8) -- (s8) -- (t8) -- (r8);

\draw (q0) [fill=white] circle (\vr);
\draw (q1) [fill=white] circle (\vr);
\draw (q2) [fill=white] circle (\vr);
\draw (q3) [fill=white] circle (\vr);
\draw (q4) [fill=white] circle (\vr);
\draw (q5) [fill=white] circle (\vr);
\draw (q6) [fill=white] circle (\vr);
\draw (q7) [fill=white] circle (\vr);
\draw (q8) [fill=white] circle (\vr);
\draw (q9) [fill=white] circle (\vr);

\draw (p1) [fill=white] circle (\vr);
\draw (p2) [fill=white] circle (\vr);
\draw (p3) [fill=white] circle (\vr);
\draw (p4) [fill=white] circle (\vr);
\draw (p5) [fill=white] circle (\vr);
\draw (p6) [fill=white] circle (\vr);
\draw (p7) [fill=white] circle (\vr);
\draw (p8) [fill=white] circle (\vr);

\draw (r1) [fill=white] circle (\vr);
\draw (r2) [fill=white] circle (\vr);
\draw (r3) [fill=white] circle (\vr);
\draw (r4) [fill=white] circle (\vr);
\draw (r5) [fill=white] circle (\vr);
\draw (r6) [fill=white] circle (\vr);
\draw (r7) [fill=white] circle (\vr);
\draw (r8) [fill=white] circle (\vr);

\draw (s1) [fill=white] circle (\vr);
\draw (s2) [fill=white] circle (\vr);
\draw (s3) [fill=white] circle (\vr);
\draw (s4) [fill=white] circle (\vr);
\draw (s5) [fill=white] circle (\vr);
\draw (s6) [fill=white] circle (\vr);
\draw (s7) [fill=white] circle (\vr);
\draw (s8) [fill=white] circle (\vr);

\draw (t1) [fill=white] circle (\vr);
\draw (t2) [fill=white] circle (\vr);
\draw (t3) [fill=white] circle (\vr);
\draw (t4) [fill=white] circle (\vr);
\draw (t5) [fill=white] circle (\vr);
\draw (t6) [fill=white] circle (\vr);
\draw (t7) [fill=white] circle (\vr);
\draw (t8) [fill=white] circle (\vr);

\draw (u1) [fill=white] circle (\vr);
\draw (u2) [fill=white] circle (\vr);
\draw (u3) [fill=white] circle (\vr);
\draw (u4) [fill=white] circle (\vr);
\draw (u5) [fill=white] circle (\vr);
\draw (u6) [fill=white] circle (\vr);
\draw (u7) [fill=white] circle (\vr);
\draw (u8) [fill=white] circle (\vr);

\draw (v1) [fill=white] circle (\vr);
\draw (v2) [fill=white] circle (\vr);
\draw (v3) [fill=white] circle (\vr);
\draw (v4) [fill=white] circle (\vr);
\draw (v5) [fill=white] circle (\vr);
\draw (v6) [fill=white] circle (\vr);
\draw (v7) [fill=white] circle (\vr);
\draw (v8) [fill=white] circle (\vr);

\draw[anchor = south] (v1) node {$7$};
\draw[anchor = south] (v2) node {$8$};
\draw[anchor = south] (v3) node {$1$};
\draw[anchor = south] (v4) node {$2$};
\draw[anchor = south] (v5) node {$3$};
\draw[anchor = south] (v6) node {$4$};
\draw[anchor = south] (v7) node {$5$};
\draw[anchor = south] (v8) node {$6$};

\draw[anchor = south] (u1) node {$6$};
\draw[anchor = south] (u2) node {$7$};
\draw[anchor = south] (u3) node {$8$};
\draw[anchor = south] (u4) node {$1$};
\draw[anchor = south] (u5) node {$2$};
\draw[anchor = south] (u6) node {$3$};
\draw[anchor = south] (u7) node {$4$};
\draw[anchor = south] (u8) node {$5$};

\draw[anchor = south] (t1) node {$5$};
\draw[anchor = south] (t2) node {$6$};
\draw[anchor = south] (t3) node {$7$};
\draw[anchor = south] (t4) node {$8$};
\draw[anchor = south] (t5) node {$1$};
\draw[anchor = south] (t6) node {$2$};
\draw[anchor = south] (t7) node {$3$};
\draw[anchor = south] (t8) node {$4$};

\draw[anchor = west] (r1) node {$4$};
\draw[anchor = west] (r2) node {$5$};
\draw[anchor = west] (r3) node {$6$};
\draw[anchor = west] (r4) node {$7$};
\draw[anchor = west] (r5) node {$8$};
\draw[anchor = west] (r6) node {$1$};
\draw[anchor = west] (r7) node {$2$};
\draw[anchor = west] (r8) node {$3$};

\draw[anchor = east] (s1) node {$4$};
\draw[anchor = east] (s2) node {$5$};
\draw[anchor = east] (s3) node {$6$};
\draw[anchor = east] (s4) node {$7$};
\draw[anchor = east] (s5) node {$8$};
\draw[anchor = east] (s6) node {$1$};
\draw[anchor = east] (s7) node {$2$};
\draw[anchor = east] (s8) node {$3$};

\draw[anchor = west] (p1) node {$3$};
\draw[anchor = west] (p2) node {$4$};
\draw[anchor = west] (p3) node {$5$};
\draw[anchor = west] (p4) node {$6$};
\draw[anchor = west] (p5) node {$7$};
\draw[anchor = west] (p6) node {$8$};
\draw[anchor = west] (p7) node {$1$};
\draw[anchor = west] (p8) node {$2$};

\draw[anchor = east] (q0) node {$8$};
\draw (1.1,-0.2) node {$1$};
\draw (2.1,-0.2) node {$2$};
\draw (3.1,-0.2) node {$3$};
\draw (4.1,-0.2) node {$4$};
\draw (5.1,-0.2) node {$5$};
\draw (6.1,-0.2) node {$6$};
\draw (7.1,-0.2) node {$7$};
\draw (8.1,-0.2) node {$8$};
\draw[anchor = west] (q9) node {$1$};

\end{tikzpicture}
\end{center}
\caption{Graph $G_2$ with $\Delta(G_2)=5<8=A_b(G_2)$.}
\label{G_5}
\end{figure}

Notice that $d_{G_n}(y_1^i)=n+3$, $d_{G_n}(x_j^i)=n+2$, $j\in[2]$, and that the degree of the other vertices of $G_n$ is at most $2$. Hence, $\Delta(G_n)=n+3$ and $m(G_n)=n+4$ (actually, we have also $\varphi(G_n)=n+4$, the appropriate b-coloring is easy to obtain). Next we show that $A_b(G_n)=2n+4$. For given $i\in[2n+4]$, observe the partition $P=\{A_0^i=\{y_1^{i-1}, y_1^{i+1}, z_j^i:j\in[n-1]\},A_1^i=\{x_1^i,x_2^i\}\}$ of $N_{G_n}(y_1^i)$,  where $y_1^0=u$ when $i=1$ and $y_1^{2n+5}=v$ when $i=2n+4$. First notice that ${\rm elp}_{G_n}(y_1^i,P)=n+1$ as $\{x_1^iy_j^ix_2^i:j\in [2,n+2]\}$ are  EVI disjoint $(x_1^i,x_2^i)$-paths of even length. Thus, $d_{G_n}^a(y_1^i)=(n+1)+1+(n+1)=2n+3$, $i\in[2n+4]$, which gives $m_a(G_n)=2n+4$, as there are exactly $2n+4$ vertices of this acyclic degree. By Theorem \ref{madegree} we have $A_b(G_n)\leq 2n+4$.

To show the equality we construct an acyclic b-coloring $c:V(G_n)\rightarrow[2n+4]$ with an acyclic b-vertex in every color class. For that purpose, let
$$A_i=\{c(x_1^i)=c(x_2^i), c(y_2^i), \dots, c(y_{n+2}^i), c(z_1^i), \dots, c(z_{n-1}^i)\}$$
and
$$
B_i=([2n+4]-[i-1,i+1])-\{-2n-4+i+1,2n+4+i-1\}
$$
for $i\in [2n+4]$.
Notice that $|A_i|=|B_i|=2n+1$. If $c$ satisfies $c(y_1^i)=i$ and $c(x_1^i)=c(x_2^i)$ for every $i\in [2n+4]$, $A_i=B_i$ for every $i\in [2n+4]$, $c(u)=2n+4$ and $c(v)=1$, then every vertex $y_1^i$, $i\in [2n+4]$, is an acyclic b-vertex. Hence, $A_b(G_n)\geq 2n+4$. One such coloring for $G_2$ is presented in Figure \ref{G_5}. With this $A_b(G_n)=2n+4$ and $A_b(G_n)-\Delta(G_n)=n+1$ for every $n\geq 1$ and finally $(A_b(G_n)-\Delta(G_n))\rightarrow\infty$ as $n\rightarrow\infty$. 
\qed\bigskip

In the remainder of this section we will prove that there is a nonlinear bound on $A_b(G)$ with respect to $\Delta(G)$ and present an infinite family of extremal graphs for this bound. The following upper bound on $A_b(G)$ can be deduced from Theorem \ref{madegree}.
\begin{corollary}\label{quadraticBound}
For any graph $G$ with $\Delta(G)\geq 2$ we have $A_b(G)\leq \frac{1}{2}(\Delta(G))^2+1$.
\end{corollary}

\noindent {\textbf{Proof.}} 
Observe that for any weak partition $P=\{A^P_0,A^P_1,\dots,A^P_k\}$ of $N_G(v)$ there is $$(|P|-1)\leq \left\lfloor\frac{d_G(v) - |A_0^p|}2\right\rfloor\;\;\;\textrm{and}\;\;\; {\rm elp}_G(v,P)\leq \left\lfloor\frac{d_G(v) - |A_0^p|}2\right\rfloor(\Delta(G) -1).
$$
This implies
$$
d_G^a(v)\leq |A^P_0|+\frac{d_G(v) - |A_0^p|}{2}+\frac{d_G(v) - |A_0^p|}{2}(\Delta(G)-1)\leq \frac{1}{2}(\Delta(G))^2.
$$
In consequence
$$
m_a(G)\leq \frac{1}{2}(\Delta(G))^2+1.
$$
\qed\bigskip

The above bound is tight, as the example of $C_4$ shows. Moreover, it belongs to an infinite family of graphs satisfying the condition with equality.

\begin{theorem}\label{DeltaSquaredExtremal}
There exists an infinite family of graphs $G_1, G_2, \dots$ such that $A_b(G_n)=m_a(G_n)=\frac{1}{2}(\Delta(G_n))^2+1$. 
\end{theorem}

\begin{figure}[ht!]
\begin{center}
\begin{tikzpicture}[scale=1,style=thick,x=1cm,y=1cm]
\def\vr{2pt} % \vr = vertex radius;

% define vertices
%%%%%
%%%%%
\path (0.5,4) coordinate (a);
\path (4.25,4) coordinate (b);
\path (10.5,4) coordinate (c);
\path (0,2) coordinate (d);
\path (1,2) coordinate (e);
\path (2.5,2) coordinate (f);
\path (3.5,2) coordinate (g);
\path (5,2) coordinate (h);
\path (6,2) coordinate (i);
\path (8,2) coordinate (j);
\path (9,2) coordinate (k);
\path (10.5,2) coordinate (l);
\path (11.5,2) coordinate (m);
\path (13,2) coordinate (n);
\path (14,2) coordinate (o);
\path (0.5,0) coordinate (p);
\path (2.25,0) coordinate (r);
\path (3,0) coordinate (s);
\path (3.75,0) coordinate (t);
\path (4.75,0) coordinate (u);
\path (5.5,0) coordinate (v);
\path (6.25,0) coordinate (w);
\path (7,0) coordinate (z);
\path (7.5,0) coordinate (a1);
\path (8,0) coordinate (b1);
\path (8.5,0) coordinate (c1);
\path (9,0) coordinate (d1);
\path (9.75,0) coordinate (e1);
\path (10.25,0) coordinate (f1);
\path (10.75,0) coordinate (g1);
\path (11.25,0) coordinate (h1);
\path (11.75,0) coordinate (i1);
\path (12.5,0) coordinate (j1);
\path (13,0) coordinate (k1);
\path (13.5,0) coordinate (l1);
\path (14,0) coordinate (m1);
\path (14.5,0) coordinate (n1);

\draw (a) -- (d) -- (p);
\draw (a) -- (e) -- (p);
\draw (b) -- (f) -- (r);
\draw (b) -- (g) -- (r);
\draw  (f) -- (s);
\draw  (g) -- (s);
\draw  (f) -- (t);
\draw  (g) -- (t);
\draw (b) -- (h) -- (v);
\draw (b) -- (i) -- (v);
\draw  (h) -- (u);
\draw  (i) -- (u);
\draw  (h) -- (w);
\draw  (i) -- (w);
\draw (c) -- (j) -- (z);
\draw (c) -- (k) -- (z);
\draw (j) -- (a1);
\draw (k) -- (a1);
\draw (j) -- (b1);
\draw (k) -- (b1);
\draw (j) -- (c1);
\draw (k) -- (c1);
\draw (j) -- (d1);
\draw (k) -- (d1);
\draw (c) -- (l) -- (e1);
\draw (c) -- (m) -- (e1);
\draw (l) -- (f1);
\draw (m) -- (f1);
\draw (l) -- (g1);
\draw (m) -- (g1);
\draw (l) -- (h1);
\draw (m) -- (h1);
\draw (l) -- (i1);
\draw (m) -- (i1);
\draw (c) -- (n) -- (j1);
\draw (c) -- (o) -- (j1);
\draw (n) -- (k1);
\draw (o) -- (k1);
\draw (n) -- (l1);
\draw (o) -- (l1);
\draw (n) -- (m1);
\draw (o) -- (m1);
\draw (n) -- (n1);
\draw (o) -- (n1);

\draw (a) [fill=white] circle (\vr);
\draw (b) [fill=white] circle (\vr);
\draw (c) [fill=white] circle (\vr);
\draw (d) [fill=white] circle (\vr);
\draw (e) [fill=white] circle (\vr);
\draw (f) [fill=white] circle (\vr);
\draw (g) [fill=white] circle (\vr);
\draw (h) [fill=white] circle (\vr);
\draw (i) [fill=white] circle (\vr);
\draw (j) [fill=white] circle (\vr);
\draw (k) [fill=white] circle (\vr);
\draw (l) [fill=white] circle (\vr);
\draw (m) [fill=white] circle (\vr);
\draw (n) [fill=white] circle (\vr);
\draw (o) [fill=white] circle (\vr);
\draw (p) [fill=white] circle (\vr);
\draw (r) [fill=white] circle (\vr);
\draw (s) [fill=white] circle (\vr);
\draw (t) [fill=white] circle (\vr);
\draw (u) [fill=white] circle (\vr);
\draw (v) [fill=white] circle (\vr);
\draw (w) [fill=white] circle (\vr);
\draw (z) [fill=white] circle (\vr);
\draw (a1) [fill=white] circle (\vr);
\draw (b1) [fill=white] circle (\vr);
\draw (c1) [fill=white] circle (\vr);
\draw (d1) [fill=white] circle (\vr);
\draw (e1) [fill=white] circle (\vr);
\draw (f1) [fill=white] circle (\vr);
\draw (g1) [fill=white] circle (\vr);
\draw (h1) [fill=white] circle (\vr);
\draw (i1) [fill=white] circle (\vr);
\draw (j1) [fill=white] circle (\vr);
\draw (k1) [fill=white] circle (\vr);
\draw (l1) [fill=white] circle (\vr);
\draw (m1) [fill=white] circle (\vr);
\draw (n1) [fill=white] circle (\vr);

\draw[anchor = south] (a) node {$v^i$};
\draw[anchor = south] (b) node {$v^i$};
\draw[anchor = south] (c) node {$v^i$};
\draw[anchor = east] (d) node {$x_1^i$};
\draw[anchor = west] (e) node {$x_2^i$};
\draw[anchor = east] (f) node {$x_1^i$};
\draw[anchor = west] (g) node {$x_2^i$};
\draw[anchor = east] (h) node {$x_3^i$};
\draw[anchor = west] (i) node {$x_4^i$};
\draw[anchor = east] (j) node {$x_1^i$};
\draw[anchor = west] (k) node {$x_2^i$};
\draw[anchor = east] (l) node {$x_3^i$};
\draw[anchor = west] (m) node {$x_4^i$};
\draw[anchor = west] (n) node {$x_5^i$};
\draw[anchor = west] (o) node {$x_6^i$};
\draw[anchor = north] (p) node {$y_1^i$};
\draw[anchor = north] (r) node {$y_1^i$};
\draw[anchor = north] (s) node {$y_2^i$};
\draw[anchor = north] (t) node {$y_3^i$};
\draw[anchor = north] (u) node {$y_4^i$};
\draw[anchor = north] (v) node {$y_5^i$};
\draw[anchor = north] (w) node {$y_6^i$};
\draw[anchor = north] (z) node {$y_1^i$};
\draw[anchor = north] (a1) node {$y_2^i$};
\draw[anchor = north] (b1) node {$y_3^i$};
\draw[anchor = north] (c1) node {$y_4^i$};
\draw[anchor = north] (d1) node {$y_5^i$};
\draw[anchor = north] (e1) node {$y_6^i$};
\draw[anchor = north] (f1) node {$y_7^i$};
\draw[anchor = north] (g1) node {$y_8^i$};
\draw[anchor = north] (h1) node {$y_9^i$};
\draw[anchor = north] (i1) node {$y_{10}^i$};
\draw[anchor = north] (j1) node {$y_{11}^i$};
\draw[anchor = north] (k1) node {$y_{12}^i$};
\draw[anchor = north] (l1) node {$y_{13}^i$};
\draw[anchor = north] (m1) node {$y_{14}^i$};
\draw[anchor = north] (n1) node {$y_{15}^i$};
\end{tikzpicture}
\end{center}
\caption{Graphs $H_{1,i}=C_4$, $H_{2,i}$ and $H_{3,i}$.}\label{higHni}
\end{figure}
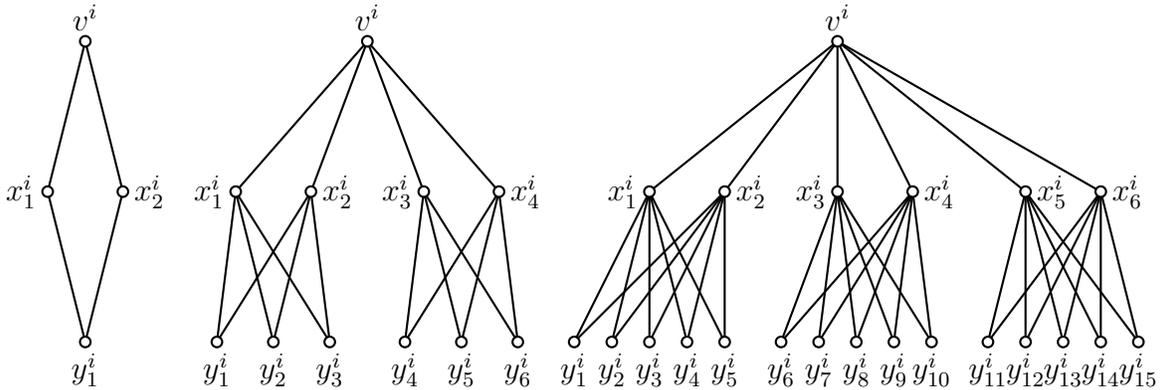

\noindent {\textbf{Proof.}}
For any positive integer $n$, we define a graph $H_n$. Then, by combining some number of copies $H_{n,i}$ of $H_n$, we will define graph $G_n$ being a member of the desired family. The vertices and edges of $H_{n,i}$ are defined as follows: 
\begin{align*}
V(H_{n,i})&=\{v^i\}\cup\{x^i_j:j\in[2n]\}\cup\{y^i_k:k\in[n(2n-1)]\},\\
E(H_{n,i})&=\{v^ix^i_j:j\in[2n]\}\\
&\cup\{x^i_jy^i_k:j\in\{2\ell-1,2\ell\}, k\in[(2n-1)(\ell-1)+1,(2n-1)\ell], \ell\in [n]\}.
\end{align*}
Graphs $H_{1,i}=C_4$, $H_{2,i}$ and $H_{3,i}$ are presented in Figure \ref{higHni}.

Now let $G_1=H_{1,1}=C_4$. Recall that $A_b(G_1)=3=\frac{1}{2}(\Delta(G_1))^2+1$ and, since for every $v\in V(G_1)$ we have $d_{G_1}^a(v)=2$, also $m_a(G_1)=3$. For $n\geq 2$, in order to obtain $G_n$ we take $2n^2+1$ graphs $H_{n,i}$, $i\in[0,2n^2]$, and identify $y_{n(2n-1)}^{i-1}$ with $y_1^{i}$ for $i\in[2n^2]$. Note that $G_n$ has $2n^2+1$ vertices $v^i$, $i\in[0,2n^2]$ of degree $2n$, $(2n^2+1)2n$ vertices $x_j^i$, $i\in[0,2n^2]$, $j\in[2n]$ of degree $2n$, $2n^2$ vertices $y_{n(2n-1)}^{i-1}=y_1^{i}$, $i\in[2n^2]$, of degree $4\leq 2n$ and $n(2n-1)(2n^2+1)-4n^2=(2n^2+1)(n(2n-1)-2)+1$ vertices $y_k^i$, $(i,k)\in([0,2n^2]\times[2,n(2n-1)-1])\cup\{(0,1),(2n^2,n(2n-1))\}$ of degree $2$. Thus $\Delta(G)=2n$.

For every vertex $v^i$, $i\in[0,2n^2]$, we have $d_{G_n}^a(v^i)=2n^2$. Indeed, for the weak partition $P$ of its neighborhood $A_0=\emptyset$, $A_{\ell}=\{x_{2\ell-1}^i,x_{2\ell}^i\}$ for $\ell\in[n]$, we have $|A_0^P|=0$, $|P|-1=n$ and ${\rm elp}_{G_n}(v,P)=n(2n-1)$ is the number of EVI disjoint paths of the form $x_{2\ell-1}^iy_k^ix_{2\ell}^i$, where $(\ell,k)\in[n]\times[(2n-1)(\ell-1)+1,(2n-1)\ell]$. This implies that $d_{G_n}^a(v^i)\geq 0+n+n(2n-1)=2n^2=\frac{1}{2}(\Delta(G_n))^2$. The inequality $d_{G_n}^a(v^i)\leq\frac{1}{2}(\Delta(G_n))^2$ follows from the proof of Corollary \ref{quadraticBound}. Obviously no other vertex can have higher acyclic degree. Since there are $2n^2+1$ vertices $v^i$, from Theorem \ref{madegree} we obtain $A_b(G_n)\leq m_a(G_n)=2n^2+1$.

In order to finish the proof we define the following coloring $c$, using the elements of additive group $\mathbb{Z}_{2n^2+1}$ as colors (with addition modulo $(2n^2+1)$ in the formulae):
\begin{align*}
    c(y^i_1)&=i, &&i\in[0,2n^2],\\
    c(y^i_{n(2n-1})&=i+1, &&i\in[0,2n^2],\\
    c(y^i_k)&=i+k, &&i\in[0,2n^2],  k\in[2,n(2n-1)-1],\\
    c(x^i_{2\ell-1})=c(x^i_{2\ell})&=i+n(2n-1)-1+\ell, &&i\in[0,2n^2],  \ell\in[n],\\
    c(v^i)&=i+2n^2, &&i\in[0,2n^2].\\
\end{align*}
Note that $c$ is an acyclic coloring of $G_n$ (only the colors of $x_j^i$ are used twice in every part $H_{n,i}$). Moreover, every $v^i$ is an acyclic b-vertex: $n$ colors are blocked by its neighbors $x_j^i$ and $n(2n-1)$ other colors because of the cycles $v^ix_{2\ell-1}^iy_k^ix_{2\ell}^iv^i$, which makes any of the potential $2n^2$ recolorings impossible. Since there is a vertex $v^i$ in every color class, $c$ is an acyclic b-coloring and $A_b(G_n)\geq 2n^2+1$. \qed\bigskip

\begin{remark}
Note that the family defined in the proof of Theorem \ref{DeltaSquaredExtremal} is another family proving Theorem \ref{roofsInfinite}.
\end{remark}

%%%%%%%%%%%%%%%%%%%%%%%%%%%%%%%%%%%%%%%%%%%%%%%%%%%%%%%
%%%%%%%%%%%%%%%%%%%%%%%%%%%%%%%%%%%%%%%%%%%%%%%%%%%%555

\section{Acyclic b-chromatic number of joins}\label{sec_BAcyclicJoins}

For more exact results we recall that the join of graphs $G$ and $H$ is the graph $G\vee H$ obtained from disjoint copies of $G$ and $H$ joined with all the possible edges between $V(G)$ and $V(H)$. More formally, $V(G\vee H)=V(G)\sqcup V(H)$ and $E(G\vee H)=E(G)\sqcup E(H)\sqcup \{uv:u\in V(G)\wedge v\in V(H)\}$ where $\sqcup $ denotes the disjoint union.

\begin{theorem}\label{join}
For two non-complete graphs $G$ and $H$ we have 
$$A_b(G\vee H)=\max\{A_b(G)+n_H,A_b(H)+n_G\}.$$
If $H\cong K_q$, then $A_b(G\vee H)=A_b(G)+q$.
\end{theorem}

\noindent {\textbf{Proof.}}
Let $G$ and $H$ be two non-complete graphs. Let $c_G$ be an $A_b(G)$-coloring of $G$ and let $V(H)=\{v_1,\dots, v_{n_H}\}$. The map 
 \begin{equation}\label{coloring}
c(v)=\left\{ 
\begin{array}{lcl}
c_G(v) & \text{if} & v\in V(G), \\
A_b(G)+i & \text{if} & v=v_{i}\in V(H), 
\end{array}
\right. 
\end{equation}
is an acyclic coloring of $G\vee H$, because $c_G$ is an acyclic coloring and all vertices of $H$ receive different colors. Suppose that we can present a recoloring step for color $i$ of $c$ in $G\vee H$ to obtain coloring $c'$. If $i>A_b(G)$, then vertex $v_{i-A_b(G)}\in V(H)$ is the only vertex of color $i$. Since $v_{i-A_b(G)}$ is adjacent to all vertices of $G$, $i$ is recolored in $c'$ by some color $j>A_b(G)$ where $j\neq i$. Again, $v_{j-A_b(G)}\in V(H)$ is the only vertex from $H$ with color $j$ in coloring $c$. Since $G$ is not complete, we have $A_b(G)<n_G$ by Proposition \ref{order} and there exist $x,y\in V(G)$ with $c(x)=c(y)$. But now $xv_{i-A_b(G)}yv_{j-A_b(G)}x$ is a bi-chromatic $4$-cycle under coloring $c'$, a contradiction. So we may assume that $i\leq A_b(G)$ and that all the vertices of color $i$ are from $V(G)$. Every vertex of color $i$ is adjacent to all vertices of $H$ and they can therefore be recolored only with colors already used in $G$. But this is not possible because $c_G$ is an $A_b(G)$-coloring. Hence $A_b(G\vee H)\geq A_b(G)+n_H$. By symmetric arguments we can show that $A_b(G\vee H)\geq A_b(H)+n_G$ which yields $A_b(G\vee H)\geq \max\{A_b(G)+n_H,A_b(H)+n_G\}$.      

We prove the opposite inequality by a contradiction. Indeed, suppose that $A_b(G\vee H) > \max\{A_b(G)+n_H, A_b(H)+n_G\}$ for some graphs $G$ and $H$ that are not complete. Let $c$ be an $A_b(G\vee H)$-coloring and let $c_G$ and $c_H$ be colorings of $G$ and $H$, respectively, induced by $c$, that is $c_G(v)=c(v)$ for every $v\in V(G)$ and $c_H(u)=c(u)$ for every $u\in V(H)$. Clearly colors of $c_G$ are different than colors of $c_H$. If $c_G$ has less than $n_G$ colors and $c_H$ less than $n_H$ colors, then $c_G(u)=c_G(v)$ and $c_H(x)=c_H(y)$ for some $u,v\in V(G)$ and $x,y\in V(H)$. But then we have a bi-colored four cycle $uxvyu$, a contradiction. So, $c_G$ has $n_G$ colors or $c_H$ has $n_H$ colors. Without loss of generality we can assume that $c_G$ has $n_G$ colors. But then $c_H$ has more than $A_b(H)$ colors in $H$ and does not yield a minimal element with respect to $\prec_a$ in $\mathcal{AF}(H)$. So, there exists $c'_H$ such that $c'_H\triangleleft_a c_H$ and the coloring 
 \begin{equation*}
c'(v)=\left\{ 
\begin{array}{lcl}
c(v) & \text{if} & v\in V(G), \\
c'_H(v) & \text{if} & v\in V(H), 
\end{array}
\right. 
\end{equation*}
is obtained from $c$ by an acyclic recoloring step, a contradiction with $c$ being an $A_b(G\vee H)$-coloring. This yields the desired equality and we are done with the first part. 

If $G\cong K_p$ and $H\cong K_q$, then $G\vee H\cong K_{p+q}$ and equality holds by Corollary \ref{order}. If only one of $G$ and $H$ is complete, say $H\cong K_q$, then we can use the coloring $c$ defined in (\ref{coloring}) for an $A_b(G)$-coloring $c_G$ of $G$. Following the same reasoning after (\ref{coloring}), this time only for $i\leq A_b(G)$, we obtain that $A_b(G\vee K_q)\geq A_b(G)+q$ (notice that $i>A_b(G)$ yields a contradiction since $H\cong K_q$ now). 

Conversely, suppose that $A_b(G\vee K_q)>A_b(G)+q$ and let $c$ be an $A_b(G\vee H)$-coloring. Again we can follow above steps and see that $c_G$, that is the restriction of $c$ to $G$, contains more than $A_b(G)$ colors and one can perform an acyclic recoloring step in $G$. This yields an acycling recoloring step in $G\vee K_q$ and $c$ is not a minimal element of $\prec_a$, a final contradiction.\qed\bigskip

Recall that complete bipartite graph $K_{m,n}=\overline{K}_n\vee\overline{K}_m$, wheel $W_n=K_1\vee C_{n-1}$, fan $F_n=K_1\vee P_{n-1}$ and complete split graph $K_n\vee\overline{K}_m$ are all joins of two graphs. Hence the following corollary follows directly from  Theorem \ref{join} and Corollaries \ref{order} and \ref{basic2}.

\begin{corollary}\label{wheel}
For every positive integers $k,\ell,m,n$, where $k,\ell\geq 5$,  we have 
\begin{itemize}
\item $A_b(K_{n,m})=1+\max\{n,m\}$;
\item $A_b(W_k)=4$;
\item $A_b(F_k)=4$;
\item $A_b(K_n\vee\overline{K}_m)=n+1$;
\item $A_b(P_k\vee P_{\ell })=A_b(P_k\vee C_{\ell })=A_b(C_k\vee C_{\ell })=3+\max\{k,\ell\}$.
\end{itemize}
\end{corollary}

Corollary \ref{wheel} implies in particular that the difference between $A_b(G)$ and $\varphi(G)$ can be arbitrarily large.
\begin{corollary}\label{wheelInfinite}
There exists an infinite family of graphs $G_1, G_2, \dots$ such that $(A_b(G_n)-\varphi(G_n))\rightarrow\infty$ as $n\rightarrow\infty$. 
\end{corollary}
\noindent {\textbf{Proof.}}
For $n\geq 1$, let $G_n=K_{n,n}$. Since the proper $2$-coloring of $K_{n,n}$ is its b-coloring, we have $\varphi(K_{n,n})=2$, while $A_b(K_{n,n})=1+n$, so $A_b(K_{n,n})-\varphi(K_{n,n})=n-1$ for $n\geq 1$ and $(A_b(G_n)-\varphi(G_n))\rightarrow\infty$ as $n\rightarrow\infty$.
\qed\bigskip

Notice also that the family $\{G_n\}$ defined in the proof of Theorem \ref{roofsInfinite} is another family of graphs for which $(A_b(G_n)-\varphi(G_n))\rightarrow\infty$ when $n\rightarrow\infty$.

%%%%%%%%%%%%%%%%%%%%%%%%%%%%%%%%%%%%%%%%%%%%%55
%%%%%%%%%%%%%%%%%%%%%%%%%%%%%%%%%%%%%%%%%%%%%%%%

\section{Final Remarks}\label{sec_FinalRemarks}

In the paper we introduced a new graph invariant, the b-acyclic chromatic number $A_b(G)$ and proved some of its properties. Some problems, however, remain still open.

Although the construction of a b-acyclic coloring seems very similar to the one of b-colorings (one considers only acyclic colorings instead of all colorings), we observed some interesting differences between them. In particular, $A_b(G)$ can be arbitrarily larger than the acyclic chromatic number $A(G)$, maximum degree $\Delta(G)$ and b-chromatic number $\varphi(G)$. The last result is consistent with the intuition that $A_b(G)$ should be not less than $\varphi(G)$, since the strictly partial ordered set of acyclic colorings is obviously a subset of the strictly partial ordered set of all proper colorings. However, we cannot neither prove nor disprove it.

On the other hand, we proved the theorem allowing to verify whether a coloring is a minimal element of the strictly partial ordered set of acyclic colorings, using the criterion of the existence of a b-acyclic vertex in every color class (where the concept of b-acyclic vertex is a natural generalization of the notion of a b-vertex). We also proved an inequality analogous to $\varphi(G)\leq m(G)$. To this end, we introduced new vertex measure, the acyclic degree $d_G^a(v)\geq d_G(v)$ and resulting graph invariant $m_a(G)\geq m(G)$ allowing to define the upper bound $A_b(G)\leq m_a(G)$. But these results still do not help to prove any relation between $A_b(G)$ and $\varphi(G)$.

\begin{figure}[ht!]
\begin{center}
\begin{tikzpicture}[scale=1.8,style=thick,x=1cm,y=1cm]
\def\vr{1.3pt} % \vr = vertex radius;

% define vertices
%%%%%
%%%%%
\path (0,0) coordinate (a);
\path (2,0) coordinate (b);
\path (4,0) coordinate (c);
\path (6,0) coordinate (d);
\path (0,3) coordinate (a1);
\path (2,3) coordinate (b1);
\path (4,3) coordinate (c1);
\path (6,3) coordinate (d1);
\path (3,4.5) coordinate (e1);
\path (3,-1.5) coordinate (e);
\path (0,1.5) coordinate (a2);
\path (1.5,1.5) coordinate (b2);
\path (4.5,1.5) coordinate (c2);
\path (6,1.5) coordinate (d2);

%  edges

\draw (a) -- (b1) -- (c) -- (d1) -- (a) -- (c1) -- (d) -- (a1) -- (c) -- (e) -- (b) -- (a1) -- (e1) -- (b1) -- (d) -- (e) -- (a) -- (a2) -- (b2) -- (c2) -- (d2);
\draw (d) -- (d2) -- (d1) -- (b) -- (b2) -- (b1);
\draw (a1) -- (a2);
\draw (c1) -- (e1) -- (d1);
\draw (b) -- (c1);
\draw (c) -- (c2) -- (c1);
\draw (e) -- (e1);

\draw[bend left] (a2) to [bend right] (d2);

\draw (a) [fill=black] circle (\vr);
\draw (b) [fill=white] circle (\vr);
\draw (c) [fill=black] circle (\vr);
\draw (d) [fill=white] circle (\vr);
\draw (a1) [fill=white] circle (\vr);
\draw (a2) [fill=white] circle (\vr);
\draw (b1) [fill=white] circle (\vr);
\draw (b2) [fill=white] circle (\vr);
\draw (c1) [fill=white] circle (\vr);
\draw (c2) [fill=white] circle (\vr);
\draw (d1) [fill=black] circle (\vr);
\draw (d2) [fill=white] circle (\vr);
\draw (e1) [fill=black] circle (\vr);
\draw (e) [fill=black] circle (\vr);

\draw[anchor = north] (a) node {$1$};
\draw[anchor = north] (b) node {$2$};
\draw[anchor = north] (c) node {$3$};
\draw[anchor = north] (d) node {$4$};
\draw[anchor = south] (a1) node {$1$};
\draw[anchor = south] (b1) node {$2$};
\draw[anchor = south] (c1) node {$3$};
\draw[anchor = south] (d1) node {$4$};
\draw[anchor = east] (a2) node {$5$};
\draw[anchor = north] (e) node {$6$};
\draw[anchor = south] (e1) node {$5$};
\draw[anchor = west] (d2) node {$6$};

\draw(1.3,1.65) node {$6 (3)$};
\draw(4.65,1.65) node {$5$};
\draw(1.3,1.35) node {$x$};
\draw(4.65,1.35) node {$y$};

\draw(0.2,1.65) node {$x'$};
\draw(5.8,1.65) node {$y'$};

\draw[anchor = west] (c1) node {$z$};
\draw[anchor = east] (b1) node {$u$};
\draw[anchor = east] (b) node {$v$};

\end{tikzpicture}
\end{center}
\caption{Graph $G$ from \cite{Tuit} for which $6=\varphi(G)\leq A_b(G)$ and the minimal colorings are incomparable.}
\label{bColAndBAcyclicCol}
\end{figure}
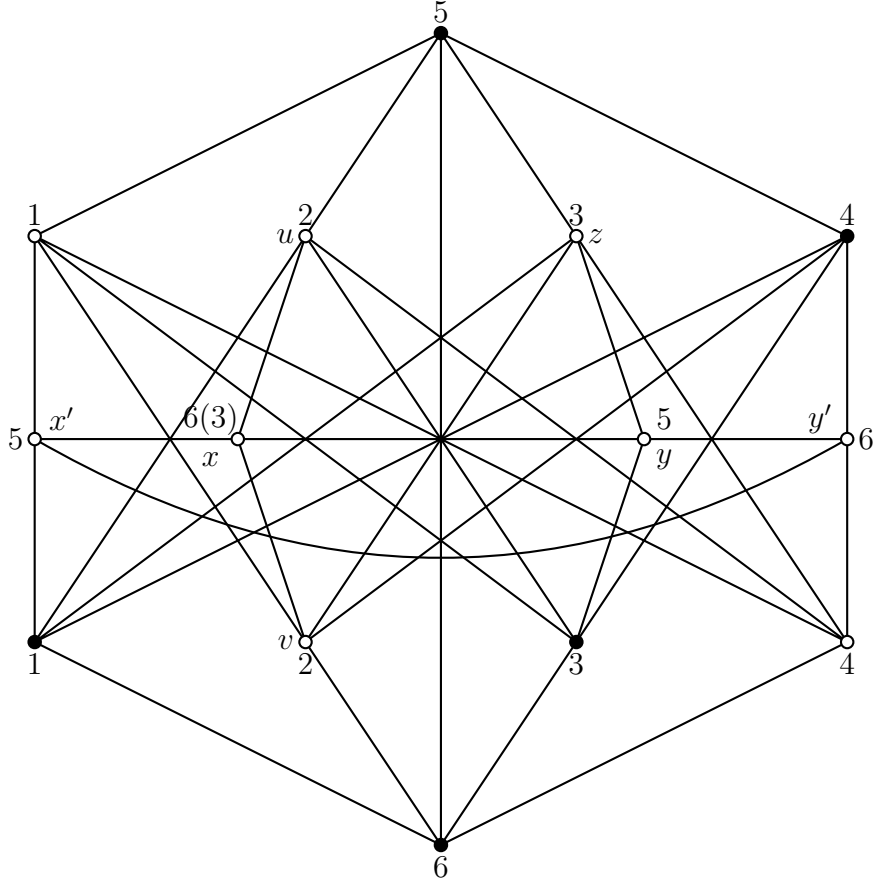

One of the causes of the difficulties in proving any relationship between these two parameters is the fact that a minimal element of the poset of proper colorings does not need to be an acyclic coloring and vice versa. A simple example is presented in Figure \ref{bColAndBAcyclicCol} (the graph was originally presented in a different context in \cite{Tuit}).

Indeed, one can easily check that the coloring $c$ in which $c(x)=6$ is the only (up to obvious rotations) b-coloring of $G$ with b-vertices $u$ and the five ones marked with black circles. However, this coloring is not acyclic, due to the cycle $xx'y'yx$. On the other hand, in the acyclic coloring $c'$ defined as $c'(x)=3$ and $c'(w)=c(w)$ for $w\in V(G)-\{x\}$, we have five b-vertices (being by definition b-acyclic vertices) marked with black circles and the sixth b-acyclic vertex $v$, that cannot be recolored with neither $1$, $3$, $4$ nor $6$ because they are the colors of its neighbors, but also color $5$ is forbidden because of the $4$-cycle $vzyxv$. On the other hand, there is no b-vertex in color class $2$. This means in particular, that $c'$ is minimal in the poset of acyclic colorings, but not in the strictly partial ordered set of proper colorings. Moreover, since $c$ and $c'$ use the same number of colors, there is no sequence of recoloring steps leading from one to the other, so they are incomparable in the poset of proper colorings.

Taking into account the above considerations, we formulate the first open problem.

\begin{problem}
Prove or disprove the inequality $A_b(G)\geq \varphi(G)$.
\end{problem}

If the inequality comes out to be false (i.e., if there is a counterexample), it would be interesting to know, for which graphs it is true. This leads us to a relaxed version of the last problem.

\begin{problem}
Characterize the graphs, for which $A_b(G)\geq \varphi(G)$.
\end{problem}

Another interesting question refers to the results about $\phi(G)$ presented in \cite{irma-99} and \cite{CaLSMaSi}.

\begin{problem}
Characterize the graphs, for which $A_b(G)\geq m_a(G)-c$ for some constant $c$. In particular, characterize the graphs, for which $A_b(G)= m_a(G)$.
\end{problem}

We have seen some examples, see Figures \ref{notenough} and \ref{system}, of $A_b(G)$-colorings that contain not recolrable critical cycle systems or, in other words, some colors have a acyclic b-vertex that is not a weak acyclic b-vertex. However all such examples also have $A_b(G)$-coloring where every color has its weak acyclic b-vertex. So, we ask if weak acyclic vertices are enough to describe acyclic b-chromatic number of a graph in the meaning of Corollary \ref{cor}?

\begin{problem}
Is it true that the acyclic b-chromatic number $A_b(G)$ of a graph $G$ is the largest integer $k$, such that there exists an acyclic $k$-coloring, where every color class $V_i$, $i\in [k]$, contains a weak acyclic b-vertex?
\end{problem}

The next problem refers to the complexity of finding $A_b(G)$, which is expected to be at lest \textbf{NP}-hard in general.

\begin{problem}
Find the complexity of deriving $A_b(G)$.
\end{problem}

We do expect that there are some special families of graphs, for which  polynomial time is enough for solving $A_b(G)$ (see Corollary \ref{tree} for trees).

\begin{problem}
Find polynomial algorithms for finding $d_G^a(v)$, $m_a(G)$ and $A_b(G)$ for chosen families of graphs.
\end{problem}

In particular some exact values or estimates for $A_b(G)$ for families of graphs other than analyzed in this paper would be a valuable contribution.

\begin{problem}
Find explicit formulae or tight lower and upper bounds on $d_G^a(v)$, $m_a(G)$ and $A_b(G)$ for chosen families of graphs.
\end{problem}

%%%%%%%%%%%%%%%%%%%%%%%%%%%%%%%%%%%%%

\end{document}